\documentclass{amsart}

\usepackage{amsmath,amssymb,amsthm,amsfonts,mathrsfs}
\usepackage[frame,cmtip,arrow,matrix,line,graph,curve]{xy}
\usepackage{graphpap,color,paralist,pstricks}
\usepackage[mathscr]{eucal}
\usepackage{mathabx}
\usepackage[pdftex,colorlinks,backref=page,citecolor=blue]{hyperref}
\usepackage{tikz}
\usepackage{epic,eepic}
\usepackage{yfonts}
\usepackage{enumerate}
\usepackage{tikz-cd}

\theoremstyle{definition}
\newtheorem{theorem}{Theorem}[section]
\newtheorem{definition}[theorem]{Definition}
\newtheorem{conjecture}[theorem]{Conjecture}
\newtheorem*{conjecturea*}{Conjecture $\operatorname{A}(X)$}
\newtheorem*{conjecturehdg*}{Conjecture $\operatorname{Hdg}(X)$}
\newtheorem{lemma}[theorem]{Lemma}
\newtheorem{proposition}[theorem]{Proposition}
\newtheorem{corollary}[theorem]{Corollary}

\newtheorem*{theorem*}{Theorem}

\theoremstyle{remark}
\newtheorem{example}[theorem]{Example}
\newtheorem{remark}[theorem]{Remark}

\numberwithin{equation}{section}

\begin{document}

\title{A decomposition theorem for Lefschetz modules}

\author{Omid Amini}
\address{\rm Omid Amini, CNRS - LMO, Université Paris-Saclay and École Polytechnique}
\email{omid.amini@universite-paris-saclay.fr}

\author{June Huh}
\address{\rm June Huh, Princeton University and Korea Institute for Advanced Study}
\email{huh@princeton.edu}

\author{Matt Larson}
\address{\rm Matt Larson, Institute for Advanced Study and Princeton University}
\email{mattlarson@princeton.edu}

\begin{abstract}
A Lefschetz module is a module over a commutative graded algebra $A$ that satisfies analogues of Poincar\'{e} duality, the hard Lefschetz property, and the Hodge--Riemann relations with respect to an open convex cone $\mathscr{K}$ in the degree one part of $A$. We analyze its decomposition into indecomposable modules over subalgebras of $A$ that are generated by elements in the closure of $\mathscr{K}$, establishing structural results that parallel the decomposition theorem for morphisms of complex projective varieties. We use our theorems to recover key statements in combinatorial Hodge theory and to illuminate the Hodge-theoretic aspects of the decomposition theorem in algebraic geometry.
\end{abstract}

\maketitle

\section{Introduction}\label{sec:introduction}

Let $A=\bigoplus_{k \ge 0} A^k$ be a  finite-dimensional commutative graded algebra over $\mathbb{R}$, and 
let $\mathscr{K}_A$ be a nonempty open convex cone in $A^1$.
Let $M= \bigoplus_{k\ge 0} M^k$ be a finite-dimensional graded $A$-module equipped with a symmetric bilinear form $\mathscr{Q} \colon M \times M \to \mathbb{R}$ that is \emph{$A$-invariant}:
\[
\mathscr{Q}(ax_1,x_2)=\mathscr{Q}(x_1,ax_2) \ \ \text{for all $a \in A$ and all $x_1,x_2 \in M$.}
\]
We say that $(M,\mathscr{Q})$ is a \emph{Lefschetz module} of degree $d$ over $(A,\mathscr{K}_A)$ if it satisfies the following three properties, called the \emph{K\"ahler package} for $(\mathscr{Q},\mathscr{K}_A)$:
\begin{enumerate}\itemsep 5pt
\item[(PD)]\label{item:PD} The induced bilinear pairing between the graded pieces
\[
\mathscr{Q} \colon M^i \times M^j \longrightarrow \mathbb{R}
\]
is nondegenerate if $i+j=d$, and it is zero if $i + j \not= d$ (\emph{Poincar\'e duality}).
\item[(HL)] \label{item:HL}For each nonnegative $k \le \frac{d}{2}$ and $\eta \in \mathscr{K}_A$, the linear map
\[
M^k \longrightarrow M^{d-k}, \qquad x \longmapsto \eta^{d-2k}  x
\]
is an isomorphism of vector spaces (\emph{hard Lefschetz property}).
\item[(HR)] \label{item:HR} For each nonnegative $k \le \frac{d}{2}$ and  $\eta \in \mathscr{K}_A$, the symmetric bilinear form
\[
M^k \times M^k \longrightarrow \mathbb{R}, \qquad (x_1,x_2) \longmapsto (-1)^k \mathscr{Q}\Big(x_1, \eta^{d-2k}  x_2\Big)
\]
is positive definite on the kernel of the linear map
\[
M^k \longrightarrow M^{d-k+1}, \qquad x \longmapsto \eta^{d-2k+1} x
\]
(\emph{Hodge--Riemann relations}).
\end{enumerate}

A classic example of a Lefschetz module of degree $d$ is given by $A = \bigoplus_{k \ge 0} H^{k,k}(X,\mathbb{R})$, where $X$ is a compact K\"ahler manifold of dimension $d$ and $H^{k,k}(X, \mathbb{R}) = H^{k,k}(X)\cap H^{2k}(X,\mathbb{R})$, viewed as a module over itself. Here
\[
\mathscr{Q}= \text{the Poincar\'e pairing on $\bigoplus_{k \ge 0} H^{k,k}(X,\mathbb{R})$} \ \ \text{and} \ \  \mathscr{K}_A=\text{the K\"ahler cone in $H^{1,1}(X,\mathbb{R})$}.
\]
In this case, the K\"ahler package for $(\mathscr{Q},\mathscr{K}_A)$ is a consequence of the Hodge theory of harmonic forms \cite[Chapter 3]{Huybrechts}.\footnote{The Hodge decomposition splits $H(X,\mathbb{C})$ into Lefschetz modules with complex coefficients.
For formulations of our main results over  $\mathbb{C}$ instead of $\mathbb{R}$, see Theorems~\ref{thmHomcomplex} and ~\ref{thmDecompositionHS}. For a discussion of our main results over  $\mathbb{Q}$ instead of $\mathbb{R}$, see Examples~\ref{ex:Q-nonuniqueness} and ~\ref{ex:Q-nonuniqueness2}.}  
In addition, assuming Grothendieck's standard conjectures on algebraic cycles \cite{Grothendieck}, the ring of algebraic cycles modulo numerical equivalence on a smooth projective variety is a Lefschetz module with respect to the ample cone (Proposition~\ref{prop:standardlefschetz}).
Further examples of Lefschetz modules include the various ``intersection cohomology without spaces'' discussed below. For a recent survey on this topic, see \cite{BP26}.

\begin{example}
\
\begin{enumerate}[(1)]\itemsep 5pt

    \item The combinatorial intersection cohomology of a convex polytope is a Lefschetz module with respect to the cone of strictly convex piecewise linear functions on its normal fan \cite{Karu}.

    \item The reduced Soergel bimodule of an element in a  Coxeter system is a Lefschetz module with respect to the open dominant cone \cite{Elias-Williamson}.

    \item The Chow ring of a matroid is a Lefschetz module with respect to the cone of strictly submodular set functions \cite{AHK}. Similarly, the augmented Chow ring of a matroid \cite{BHMPW} and the conormal Chow ring of a matroid \cite{ADH1} admit the structure of Lefschetz modules.

    \item The intersection cohomology of a matroid is a Lefschetz module with respect to the cone of positive functions \cite{BHMPW2}.

    \item The Chow ring of a unimodular quasi-projective quasilinear fan is a Lefschetz module with respect to the cone of strictly convex piecewise linear functions \cite{AP23}.
\end{enumerate}
\noindent See Example~\ref{exReverseKT} for a Lefschetz module that cannot arise from K\"ahler geometry or algebraic geometry.
\end{example}

In the remainder of this paper, we suppose that $(M,\mathscr{Q})$ is a Lefschetz module of degree $d$ over $(A,\mathscr{K}_A)$ and deduce a number of structural results.
Let $B$ be an $\mathbb{R}$-subalgebra of $A$   generated by a subset of the closure $\overline{\mathscr{K}}_A \subseteq A^1$ of $\mathscr{K}_A$, and set 
\[
\mathscr{K}_B\coloneq \text{the nonempty open convex cone given by the interior of $\overline{\mathscr{K}}_A \cap B^1$ in $B^1$.}
\]
When $\mathscr{K}_B$ contains an element that satisfies the conclusion of the hard Lefschetz property for $M$, for example when $\mathscr{K}_B$ contains an element of $\mathscr{K}_A$, we say that the graded $B$-module structure on $M$ is \emph{semismall}. When every element of $\mathscr{K}_B$ acts as zero on $M$, for example when $\mathscr{K}_B=0$, we say that the graded $B$-module structure on $M$ is \emph{trivial}.

We work in the category of finite-dimensional graded $B$-modules and degree-preserving homomorphisms of graded $B$-modules.
For a graded $B$-module $N$, we define the shifted $B$-module $N[-k]$
as the direct sum of the vector spaces $N^{i-k}$ placed in degree $i$.

 Our main results show that an arbitrary  decomposition of $M$ into indecomposable graded $B$-modules has a number of remarkable properties 
that we call the \emph{decomposition package}.  
Choose any decomposition 
\begin{equation}\label{eq:decomposition-indecomposable}
M\simeq \bigoplus_\alpha \bigoplus_{k} N_\alpha[-k]^{\oplus m(\alpha,k)},
\end{equation}
where $N_\alpha=N^0_\alpha \oplus \cdots \oplus N^{d(\alpha)}_\alpha$ are indecomposable graded $B$-modules satisfying 
\[
N^0_\alpha \neq 0, \ \ N^{d(\alpha)}_\alpha \neq 0, \ \ \text{and} \ \  \text{$N_\alpha \nsimeq N_\beta$ as graded $B$-modules for all $\alpha \neq \beta$}.  
\]
By the Krull--Schmidt theorem applied to the category of finite-dimensional graded $B$-modules \cite{Atiyah}, the collection $\{N_\alpha\}$ and the multiplicities $m(\alpha,k)$ are independent of the choice of decomposition of $M$. We then have the following properties.

\begin{theorem}[Simplicity]\label{thmHom}
For each $\alpha$, we have an isomorphism of $\mathbb{R}$-algebras
\[
\operatorname{Hom}_B(N_\alpha,N_\alpha) \simeq \mathbb{R},\ \mathbb{C},\ \text{or}\ \mathbb{H}.
\]
If $\alpha=\beta$ and $k$ is positive, or if $\alpha\neq\beta$ and
$d(\alpha)\le d(\beta)+2k$, then 
\[
\operatorname{Hom}_B(N_\alpha,N_\beta[-k])=0.
\]
\end{theorem}

Each of the three finite-dimensional division algebras over $\mathbb{R}$ arises as the endomorphism ring of some $N_{\alpha}$; see Examples~\ref{ex:endC}  and \ref{ex:endH}. 
The condition $d(\alpha)\le d(\beta)+2k$ says that the middle degree of $N_\alpha$ is at most the middle degree of $N_\beta[-k]$. 

\begin{theorem}[Duality]\label{thmPD}
For each $\alpha$, up to a nonzero real multiple, there is a unique nonzero real-valued $B$-invariant symmetric bilinear form $\mathscr{Q}_\alpha$ on $N_\alpha$ satisfying the degree orthogonality relations
\[
\mathscr{Q}_\alpha ( N^i_\alpha, N_\alpha^j)=0 \  \ \text{unless $i+j =d(\alpha)$.}
\]
This  $B$-invariant symmetric bilinear form $\mathscr{Q}_\alpha$ is nondegenerate.
\end{theorem}

We show that, for each $\alpha$, the bilinear form $\mathscr{Q}_\alpha$ gives $N_\alpha$ the structure of a Lefschetz module of degree $d(\alpha)$ over $(B,\mathscr{K}_B)$, up to a sign. 

\begin{theorem}[Hard Lefschetz]\label{thmHL}
For each nonnegative $k \le \frac{d(\alpha)}{2}$ and $\ell \in \mathscr{K}_B$, the linear map
\[
N_\alpha^k \longrightarrow N_\alpha^{d(\alpha)-k}, \qquad x \longmapsto \ell^{d(\alpha)-2k} x
\]
is an isomorphism of vector spaces.
\end{theorem}

\begin{theorem}[Hodge--Riemann]\label{thmHR}
There is a unique sign $\epsilon_\alpha=\pm 1$ such that,  for each nonnegative $k \le \frac{d(\alpha)}{2}$ and $\ell \in \mathscr{K}_B$, the symmetric bilinear form
\[
N_\alpha^k \times N_\alpha^k \longrightarrow \mathbb{R}, \qquad (x_1,x_2) \longmapsto (-1)^k \epsilon_\alpha\mathscr{Q}_\alpha\Big(x_1, \ell^{d(\alpha)-2k} x_2\Big)
\]
is positive definite on the kernel of the linear map
\[
N_\alpha^k \longrightarrow N_\alpha^{d(\alpha)-k+1}, \qquad x \longmapsto \ell^{d(\alpha)-2k+1} x.
\]
\end{theorem}

Throughout the paper, an important role will be played by an increasing filtration
\[
P=\Big( 0 \subseteq P_0 \subseteq P_1 \subseteq \cdots \subseteq P_{2d}=M\Big),
\]
called the \emph{perverse filtration} of $M$ over $B$, defined as follows:
Choose any element $\ell \in \mathscr{K}_B$, and choose any  decomposition of $M$ into a direct sum of cyclic graded $\mathbb{R}[\ell]$-modules of the form 
\[
 \bigoplus_{e=0}^{d(x)} \operatorname{span}(\ell^e x),
\]
where $x$ is an element of degree $k(x)$ in $M$. 
We then define $P_j$ as the subspace of $M$ spanned by those summands with $d(x)+2k(x) \le j$. Equivalently, we can define the perverse filtration as
\begin{equation}\label{eq:perverse-filtration-annihilators}
P_j \cap M^k = \sum_{c=0}^k \operatorname{im}\Big(M^c \xrightarrow{\ell^{k-c}} M^k\Big) \cap \operatorname{ker}\Big(M^k \xrightarrow{\ell^{j-c+1-k}} M^{j-c+1}\Big),
\end{equation}
where the kernels of $\ell^e$ are set to be zero for $e\le 0$. Informally speaking, the perverse filtration is defined so that each $P_j/P_{j-1}$ satisfies the hard Lefschetz property with respect to $\ell$.

It follows from  Theorem~\ref{thmHL} that the perverse filtration $P$ is independent of the choice of $\ell \in \mathscr{K}_B$. 
Indeed, we can refine the given decomposition \eqref{eq:decomposition-indecomposable}  to a decomposition into indecomposable graded $\mathbb{R}[\ell]$-modules and deduce from Theorem~\ref{thmHL}  that, for each $j$, we have
\begin{equation}\label{eq:perverse-filtration}
P_j=\bigoplus_{d(\alpha)+2k \le j} N_\alpha[-k]^{\oplus m(\alpha, k)}.
\end{equation}
In particular, the right-hand side does not depend on the choice of the decomposition \eqref{eq:decomposition-indecomposable}, and the left-hand side  is a graded $B$-submodule of $M$: 
\[
bP_j \subseteq P_j \ \  \text{for all $j$ and all $b\in B$.}
\]
Let $\operatorname{Gr} \coloneq \bigoplus_j P_j/P_{j-1}$ be the associated graded of the perverse filtration of $M$ over $B$. 
Note that
\[
\operatorname{Gr}=\bigoplus_{i,j}
\operatorname{Gr}^{i,j}, \ \ \text{where} \ \ 
\operatorname{Gr}^{i,j}\coloneq (M^i \cap P_j) /(M^i \cap P_{j-1}),
\]
and that $\operatorname{Gr}$ has the structure of a graded $B$-module such that
$
B^k \operatorname{Gr}^{i,j} \subseteq \operatorname{Gr}^{i+k,j}$.
By construction, $M$ and $\operatorname{Gr}$ are isomorphic as graded $B$-modules.
We consider the graded vector space 
\[
V_\alpha\coloneq \bigoplus_k V_\alpha^k, \ \ \text{where} \ \  V_\alpha^k\coloneq \operatorname{Hom}_B(N_\alpha[-k],\operatorname{Gr}^{\bullet,d(\alpha)+2k}).
\]
Let $\mathbb{D}_{\alpha}$ be the division algebra $\operatorname{Hom}_B(N_{\alpha}, N_{\alpha})$ from Theorem \ref{thmHom}, and note that $N_{\alpha}$ and $V_{\alpha}$ are modules over $\mathbb{D}_{\alpha}$. 
It follows from the vanishing statement in Theorem \ref{thmHom} that, for each $j$, $k$, and $\alpha$ satisfying $d(\alpha)+2k=j$, the natural map
$
V_\alpha^k \otimes_{\,\mathbb{D}_\alpha} N_\alpha  \rightarrow \operatorname{Gr}^{\bullet,j} 
$
is injective, and the images of these maps give a decomposition of $\operatorname{Gr}^{\bullet,j}$. 

\begin{corollary}
For each $j$, there is a canonical decomposition of $B$-modules
\[
\operatorname{Gr}^{\bullet,j}=\bigoplus_{d(\alpha)+2k=j}V_\alpha^k \otimes_{\,\mathbb{D}_\alpha} N_\alpha, \ \  \text{with  $\dim_{\mathbb{D}_\alpha} V^k_{\alpha} = m(\alpha, k)$}.
\]
In other words, $\operatorname{Gr}^{\bullet,j}$ has a unique isotypic decomposition, where for each $\alpha$, the corresponding summand is the sum of all the $B$-submodules of $\operatorname{Gr}^{\bullet,j}$ that are isomorphic to $N_\alpha[-k]$. 
\end{corollary}

In particular, in the ``semismall'' case where all the graded pieces $\operatorname{Gr}^{\bullet,j}$ are zero except for $\operatorname{Gr}^{\bullet,d}$, we have a unique isotypic decomposition of $M$ into graded $B$-modules.
This is the classical Lefschetz decomposition of $M$ when 
 $B$ is generated by a single element in $\mathscr{K}_A$.

We now formulate analogues of the relative hard Lefschetz theorem and the relative Hodge--Riemann relations using $\operatorname{Gr}$. Inspecting
\eqref{eq:perverse-filtration-annihilators}, we see that  
\[
aP_j \subseteq P_{j+2k} \ \  \text{for any $a \in A^k$.} 
\]
We therefore have an $A$-module structure $*\colon A \times \operatorname{Gr} \to \operatorname{Gr}$ such that
\[
A^k  * \operatorname{Gr}^{i,j} \subseteq \operatorname{Gr}^{i+k,j+2k} \ \  \textrm{for all $k$ and $i,j$.}
\]
Note that, when the graded $B$-module structure on $M$ is semismall, this $A$-module structure on $\operatorname{Gr}$ is trivial. 
At the other extreme, when the graded $B$-module structure on $M$ is trivial, we have $\operatorname{Gr}^{k,2k}=M^k$, and the $A$-module structure on $\operatorname{Gr}$ is that of $M$.
In any case, the symmetric bilinear form $\mathscr{Q}$ on $M$ induces a  bilinear form
\[
\underline{\mathscr{Q}} \colon \operatorname{Gr}^{i,j} \times \operatorname{Gr}^{d-i,2d-j} \longrightarrow \mathbb{R} 
\]
satisfying the degree orthogonality relations; 
see Lemmas~\ref{lem:Amod} and ~\ref{lem:formdefined}. 
Note that $\underline{\mathscr{Q}}$ is both $A$-invariant and $B$-invariant:
\[
\underline{\mathscr{Q}}(a * x_1,x_2)=
\underline{\mathscr{Q}}(x_1,a * x_2)
\ \ \text{and} \ \
\underline{\mathscr{Q}}(bx_1,x_2)=
\underline{\mathscr{Q}}(x_1,bx_2) \ \ \text{for all $a\in A$, $b\in B$, and $x_1,x_2 \in \operatorname{Gr}$}.
\]
This $\underline{\mathscr{Q}}$ is nondegenerate, giving a relative version of Poincar\'e duality; see Proposition~\ref{prop:nondegenerate}. Let $\mathscr{K}_{A/B}$ be the open convex cone $\mathscr{K}_A+B^1$ in $A^1$.

\begin{theorem}[Relative Hard Lefschetz]\label{thmRHL}
For each nonnegative $j\le d$ and $\eta \in \mathscr{K}_{A/B}$, the linear map
\[
\operatorname{Gr}^{\bullet,j} \longrightarrow \operatorname{Gr}^{\bullet+d-j,2d-j}, \qquad x \longmapsto \eta^{d-j} *  x 
\]
is an isomorphism of vector spaces.
\end{theorem}

Note that, for each $\eta \in A^1$, post-composition by $\eta*-$ defines a linear map $V_\alpha^k \longrightarrow V_\alpha^{k+1}$.

\begin{corollary}\label{cor:multiplicities}
For each nonnegative $k \le (d - d(\alpha))/2$ and $\eta \in \mathscr{K}_{A/B}$, the linear map
\[
V_\alpha^k \longrightarrow
V_\alpha^{d-d(\alpha)-k}, \qquad \varphi \longmapsto \eta^{d-d(\alpha)-2k}* \varphi
\]
is an isomorphism of vector spaces. 
In particular, for each $\alpha$, the multiplicities $m(\alpha,k)=\dim_{\mathbb{D}_\alpha} V^k_\alpha$ form a symmetric and unimodal sequence in $k$.
\end{corollary}

For each $\eta \in \mathscr{K}_{A/B}$, $\ell\in \mathscr{K}_B$, and nonnegative integers
$j \le d$ and $i \le j/2$, set
\[
\operatorname{Prim}^{i,j}_{\eta,\ell} \coloneq
\operatorname{ker}(\operatorname{Gr}^{i,j} \xrightarrow{\eta^{d-j+1}*} \operatorname{Gr}^{i+d-j+1,2d-j+2}) \cap \, \operatorname{ker}( \operatorname{Gr}^{i,j} \xrightarrow{\ell^{j-2i+1}} \operatorname{Gr}^{j-i+1,j}). 
\]
It follows from Theorem \ref{thmHL} and Theorem~\ref{thmRHL} that there is a direct sum decomposition
\begin{equation}\label{eq:primitive-decomposition}
\operatorname{Gr}=\bigoplus_{\substack{j \le d\\ i \le j/2}} \,\, \bigoplus_{\substack{0 \le s \le d-j\\ 0 \le t \le j-2i}} \eta^s *
\Bigg(\ell^t \operatorname{Prim}^{i,j}_{\eta,\ell}\Bigg).
\end{equation}

\begin{theorem}[Relative Hodge--Riemann relations]\label{thmRHR}
For each nonnegative $j \le d$ and $i \le j/2$,
and for each $\eta \in \mathscr{K}_{A/B}$ and $\ell \in \mathscr{K}_B$, 
the symmetric bilinear form on $\operatorname{Gr}^{i,j}$ defined by
\[
(x_1,x_2) \longmapsto (-1)^i \underline{\mathscr{Q}}(x_1, \eta^{d-j} * \ell^{j-2i} x_2)
\]
is positive definite when restricted to the subspace $\operatorname{Prim}^{i,j}_{\eta,\ell}$.
\end{theorem}

\begin{corollary}\label{cor:kernellefschetz}
 For each nonnegative $j \le d$ and $\eta \in \mathscr{K}_{A/B}$, the $B$-module
\[
\operatorname{Prim}_\eta^{\bullet,j} \coloneq
\operatorname{ker}( \operatorname{Gr}^{\bullet,j} \xrightarrow{\eta^{d-j+1}*} \operatorname{Gr}^{\bullet+d-j+1,2d-j+2}),
\]
equipped with the form
$\underline{\mathscr{Q}}(x_1,\eta^{d-j} * x_2)$,
is a Lefschetz module of degree $j$ over $(B,\mathscr{K}_B)$.
\end{corollary}

\begin{corollary}\label{cor:kernellefschetz2}
For each nonnegative $k \le d$ and $\ell \in \mathscr{K}_B$, the $A$-module
\[
\operatorname{Prim}_\ell^{\bullet,[k]} \coloneq
\operatorname{ker}( \operatorname{Gr}^{\bullet,2\bullet+k} \xrightarrow{\ell^{k+1}} \operatorname{Gr}^{\bullet+k+1,2\bullet+k}),
\]
equipped with the form 
$\underline{\mathscr{Q}}(x_1,\ell^{k}x_2)$,
is a Lefschetz module of degree $d-k$ over $(A,\mathscr{K}_{A/B})$.
\end{corollary}

Let $C$ be any graded $B$-subalgebra of $A$ that preserves the perverse filtration of $M$ over $B$, so that $\operatorname{Gr}$ has the structure of a graded $C$-module such that
$
C^k \operatorname{Gr}^{i,j} \subseteq \operatorname{Gr}^{i+k,j}$. Note that, unlike $B$, the graded algebra $C$ need not be generated by its degree one elements. For example, when the graded $B$-module structure on $M$ is semismall, we may take $C$ to be any graded $B$-subalgebra of $A$. We assume that $C^1$ is generated by a subset of $\overline{\mathscr{K}}_A$. Set
\[
\mathscr{K}_C\coloneq \text{the nonempty convex cone given by the interior of $\overline{\mathscr{K}}_A \cap C^1$ in $C^1$.}
\]
Because we assumed that $C^1$ is generated by a subset of $\overline{\mathscr{K}}_A$, $\mathscr{K}_C$ is open. 
Choose any decomposition of $M$ as a direct sum of indecomposable graded $C$-modules
\begin{equation}\label{eq:decomposition-indecomposable2}
M\simeq \bigoplus_\alpha \bigoplus_{k} L_\alpha[-k]^{\oplus \mu(\alpha,k)},
\end{equation}
where $L_\alpha=L^0_\alpha \oplus \cdots \oplus L^{e(\alpha)}_\alpha$ are indecomposable graded $C$-modules satisfying 
\[
L^0_\alpha \neq 0, \ \ L^{e(\alpha)}_\alpha \neq 0, \ \ \text{and} \ \  \text{$L_\alpha \nsimeq L_\beta$ as graded $C$-modules for all $\alpha \neq \beta$}.  
\]
Theorems~\ref{thmHom}, ~\ref{thmPD}, \ref{thmHL}, and \ref{thmHR} can be strengthened as follows.

\begin{theorem}\label{thmDecomposition}
For each $\alpha$, there is an isomorphism of $\mathbb{R}$-algebras
\[
\operatorname{Hom}_C(L_\alpha,L_\alpha) \simeq \mathbb{R},\ \mathbb{C},\ \text{or}\ \mathbb{H}.
\]
If $\alpha=\beta$ and $k$ is positive, or if $\alpha\neq\beta$ and
$e(\alpha)\le e(\beta)+2k$, then 
\[
\operatorname{Hom}_C(L_\alpha,L_\beta[-k])=0.
\]
Up to a nonzero real multiple, there is a unique nonzero real-valued $C$-invariant symmetric bilinear form $\mathscr{Q}_\alpha$ on $L_\alpha$ satisfying the degree orthogonality relations
\[
\mathscr{Q}_\alpha ( L^i_\alpha, L_\alpha^j)=0 \  \ \text{unless $i+j =e(\alpha)$,}
\]
and there is a unique $\epsilon_\alpha=\pm 1$ such that $(L_\alpha,\epsilon_\alpha\mathscr{Q}_\alpha)$ is a Lefschetz module of degree $e(\alpha)$ over $(C,\mathscr{K}_C)$.
\end{theorem}

The proof of Theorem~\ref{thmDecomposition} produces a distinguished isomorphism of graded $C$-modules $M \simeq \operatorname{Gr}$, whose construction depends on the choice of $\eta \in \mathscr{K}_{A/B}$.\footnote{Following ideas of Deligne \cite{Deligen1994}, in \cite{dCCanonical}, de Cataldo constructs several distinguished isomorphisms using the choice of $\eta$ in the geometric setting.}

In Section~\ref{sec:example}, we illustrate some notable features of the main theorems through explicit combinatorial examples. Examples~\ref{Fano1} and~\ref{Fano2} give a detailed analysis of a Lefschetz module obtained from the Fano plane.
 Example~\ref{VariousC} compares the conclusions of Theorem~\ref{thmDecomposition} for various possible choices of $C$.  
 Example~\ref{ExampleHirzeburch} demonstrates that the cones $\mathscr{K}$ arising from geometry are often smaller than the maximum possible.
 Example~\ref{exReverseKT} provides a Lefschetz module that cannot appear inside the cohomology of a compact K\"ahler manifold or the Chow ring of a smooth projective variety over an algebraically closed field modulo numerical equivalence. Examples~\ref{ex:endC} and~\ref{ex:endH} show that both the field of complex numbers $\mathbb{C}$ and the quaternions $\mathbb{H}$ can occur in Theorem~\ref{thmHom}. Examples~\ref{ex:Q-nonuniqueness} and ~\ref{ex:Q-nonuniqueness2} show that the uniqueness of $\mathscr{Q}_\alpha$ in Theorems~\ref{thmPD} and ~\ref{thmDecomposition} can fail if one works with coefficients in $\mathbb{Q}$ rather than in $\mathbb{R}$, and that not every such $\mathscr{Q}_\alpha$ satisfies the conclusions of Theorems~\ref{thmHR} and ~\ref{thmDecomposition}.

In Section~\ref{sec:applications}, we give a number of applications of the main theorems. Motivated by combinatorial applications, many authors have shown that certain concrete modules are Lefschetz  \cite{McMullen,Karu,Elias-Williamson,AHK,KaruRHL,BHMPW,ADH1,BHMPW2}. 
Theorems~\ref{thmHL} and ~\ref{thmHR} 
provide a tool for producing new Lefschetz modules from a given one, and some of the results cited above admit notably  simpler proofs using our work.
As an example, in Section~\ref{sec:matroids}, we give a proof of the top-heavy conjecture of Dowling--Wilson \cite{DW1,DW2} that avoids the long inductive argument of \cite{BHMPW2}.
In Section~\ref{ssec:polytopes}, as another example, we give a combinatorial proof of Stanley's result on the unimodality of the local $h$-vector \cite{subdivisions} 
that avoids Karu's relative hard Lefschetz theorem for polytopal subdivisions \cite{KaruRHL}.

Our results can also be used to study Lefschetz modules arising from algebraic geometry. In Section~\ref{sec:algebraiccycles}, we use our main theorems to verify new cases of Grothendieck's standard conjectures on algebraic cycles  \cite{Grothendieck}. In the forthcoming work \cite{LStapledon}, the authors use this to prove Fulton--Lazarsfeld positivity for matroids, settling a conjecture of Speyer and Berget--Fink. In Section~\ref{sec:hodge-decomposition-theorem}, by combining our results with \cite{BBD}, we show that the intersection cohomology of a complex projective variety admits a polarized pure Hodge structure, recovering a result of Saito \cite{Saito88, Saito90}.

The proofs of our main theorems in Sections~\ref{sec:perverse}, ~\ref{sec:relative}, and ~\ref{sec:proof} are inspired by the detailed analysis of the intersection cohomology of matroids  in \cite{BHMPW2} and by the beautiful work of de Cataldo and Migliorini on the decomposition theorem \cite{dCM02,dCM05,dCMHodge}. Although our inductive strategy resembles theirs, important differences arise. For instance, the argument of \cite{dCM05} relies on the Lefschetz hyperplane theorem and the simplicity of certain perverse sheaves,  neither of which is available in our setting.

\subsection*{Acknowledgements}
We thank Dave Anderson, Mark de Cataldo, Eduardo Cattani, Ben Elias, Leonardo Mihalcea, and Geordie Williamson for helpful conversations and insightful comments. We are especially grateful to  Christian Schnell for discussions related to Example~\ref{ExampleHirzeburch}. ChatGPT and Claude were used for proofreading. 
Part of this work was carried out at the Korea Institute for Advanced Study, and we thank the institute for a pleasant working environment.  
June Huh was partially supported by the Simons Investigator Grant.

\section{Examples}\label{sec:example}

We illustrate the decomposition package explicitly using combinatorial examples.

\begin{example}\label{Fano1}
The \emph{Fano plane} $F_7$ is the configuration of seven points and seven lines shown in the following figure:
\begin{center}
\begin{tikzpicture}[scale=1.4]
\draw (0,0) -- (2,0);
\draw (2,0) -- (1,1.73);
\draw (1,1.73) -- (0,0);
\draw (0,0) -- (1.5,.866);
\draw (2,0) -- (.5,.866);
\draw (1,1.73) -- (1,0);
\draw (1,0.577) circle [radius=0.577];

\draw [fill=black] (0,0) circle (1.5pt);      
\draw [fill=black] (2,0) circle (1.5pt);      
\draw [fill=black] (1,1.73) circle (1.5pt);   
\draw [fill=black] (1.5,.866) circle (1.5pt); 
\draw [fill=black] (.5,.866) circle (1.5pt);  
\draw [fill=black] (1,0) circle (1.5pt);      
\draw [fill=black] (1,0.57) circle (1.5pt);   

\node[below left] at (0,0) {3};
\node[below right] at (2,0) {5};
\node[above] at (1,1.73) {1};
\node[right] at (1.6,.866) {6};
\node[left] at (.4,.866) {2};
\node[below] at (1,0) {4};
\node[right] at (1.15,0.57) {7};
\end{tikzpicture}
\end{center}
The \emph{graded M\"obius algebra} of $F_7$ is the graded vector space $A=A^0 \oplus A^1 \oplus A^2 \oplus A^3$, where 
\begin{alignat*}{15}
A^0&\ =\ \mathbb{R} y_\varnothing, \\
A^1&\ =\ \mathbb{R} y_1 &&\oplus\mathbb{R} y_2 &&\oplus\mathbb{R} y_3 &&\oplus\mathbb{R} y_4 &&\oplus\mathbb{R} y_5 &&\oplus\mathbb{R} y_6 &&\oplus\mathbb{R} y_7,\\
A^2&\ =\ \mathbb{R} y_{123} &&\oplus\mathbb{R} y_{147} &&\oplus\mathbb{R} y_{156} &&\oplus\mathbb{R} y_{257} &&\oplus\mathbb{R} y_{345} &&\oplus\mathbb{R} y_{367} &&\oplus\mathbb{R} y_{246},\\
A^3&\ =\ \mathbb{R} y_E.
\end{alignat*}
The degree one basis elements $y_i$ are labeled by the seven points of $F_7$ and the degree two basis elements $y_{L}$ are labeled by the seven lines of $F_7$. The multiplicative structure is given by the rules
\[
y_\varnothing =1, \quad
y_i \cdot y_j=\begin{cases} y_{\overline{ij}}& \text{if $i\neq j$, } \\ 0& \text{if $i=j$,} \end{cases}
\ \ \text{and} \ \ 
y_i \cdot y_L=\begin{cases} y_E& \text{if $i \notin L$,} \\ 0 & \text{if $i \in L$,}\end{cases}
\]
where $\overline{ij}$ is the unique line of $F_7$ containing $i$ and $j$. 
Let  $\mathscr{K}_A$ be the set of positive linear combinations of the generators $y_1,\ldots,y_7$ in $A^1$.
Let $\mathscr{Q}\colon A \times A \to \mathbb{R}$ be 
the unique $A$-invariant symmetric bilinear pairing satisfying $\mathscr{Q}(1,y_E)=1$ and the degree orthogonality relations
\[
\mathscr{Q}(A^i,A^j)=0 \ \ \text{unless $i+j=3$.}
\]
One can then check that $(A,\mathscr{Q})$ is a Lefschetz module of degree $3$ over $(A, \mathscr{K}_A)$.
This is a special case of the statement that the intersection cohomology module of a matroid is a Lefschetz module over the graded M\"obius algebra  \cite[Theorem 1.6]{BHMPW2}.

Let $B$ be the subalgebra of $A$ generated by the elements $y_1,y_3,y_5,y_7$.\footnote{The same structure arises from any other choice of four points of $F_7$ in general position.} There are infinitely many different decompositions of $M \coloneq A$ into indecomposable graded $B$-modules, each of which has three summands. One of the three summands is common to all the decompositions: the one-dimensional summand spanned by
\[
\xi\coloneq y_1-y_2+y_3-y_4+y_5-y_6+y_7.
\]
For example, $M=S_1 \oplus S_2 \oplus S_3$, where $S_1,S_2,S_3$ are graded $B$-submodules of $M$ given by 
\begin{align*}
&S^3_1= \operatorname{span}(y_E), &S^3_2&=0,&S_3^3&=0, \\
&S^2_1=\operatorname{span}(y_{123},y_{147},y_{156},y_{257},y_{345},y_{367}), &S^2_2&= 0,&S_3^2&=\operatorname{span}(y_2 \cdot \xi),\\
&S^1_1=\operatorname{span}(y_1,y_3,y_5,y_7,y_2,y_4-y_6), &S^1_2&=\operatorname{span}(\xi),&S_3^1&=0,\\
&S^0_1=\operatorname{span}(y_{\varnothing}), &S_2^0&=0, &S_3^0&=0.
\end{align*}
We also have $M=T_1 \oplus T_2 \oplus T_3$, where $T_1,T_2,T_3$ are graded $B$-submodules of $M$ given by 
\begin{align*}
&T^3_1= \operatorname{span}(y_E), &T^3_2&=0,&T_3^3&=0, \\
&T^2_1=\operatorname{span}(y_{123},y_{147},y_{156},y_{257},y_{345},y_{367}), &T^2_2&= 0,&T_3^2&=\operatorname{span}(y_4 \cdot \xi),\\
&T^1_1=\operatorname{span}(y_1,y_3,y_5,y_7,y_4,y_2-y_6), &T^1_2&=\operatorname{span}(\xi),&T_3^1&=0,\\
&T^0_1=\operatorname{span}(y_{\varnothing}), &T_2^0&=0, &T_3^0&=0.
\end{align*}
Note that $S_1\simeq T_1$ and $S_3 \simeq T_3$ as graded $B$-modules, but $S_1 \neq T_1$ and $S_3 \neq T_3$. Note also that  $S_2=T_2$ and $S_1+S_2=T_1+T_2$. 
These identities illustrate the description (\ref{eq:perverse-filtration}) of 
the perverse filtration of $M$ over $B$, which follows from Theorem~\ref{thmHL}:
\[
P_0=P_1=0, \ \  P_2=S_2=T_2, \ \ P_3=S_1+S_2=T_1+T_2, \ \ P_4=P_5=P_6=M.
\]
One implication of the main results (Theorems~\ref{thmPD}, ~\ref{thmHL}, and ~\ref{thmHR}) is that, for any  decomposition of $M$ into indecomposable graded $B$-modules
\[
M=N_1 \oplus N_2 \oplus N_3 \ \ \text{with} \ \  \dim N_2=\dim N_3=1,
\]
the restriction of the Poincar\'e pairing $\mathscr{Q}$ on $M$ to $N_1$ is nondegenerate, and $(N_1,\mathscr{Q})$ satisfies the hard Lefschetz property and the Hodge--Riemann relations with respect to $(B,\mathscr{K}_B)$. 
These assertions can be checked, for example, for $S_1$ and $T_1$. 
Another consequence is the relative hard Lefschetz property of Theorem~\ref{thmRHL}, which says that, for any real numbers $a_1,a_3,a_5,a_7$ and any positive real numbers $a_2,a_4,a_6$, we have
\[
(a_1y_1+a_2y_2+a_3y_3+a_4y_4+a_5y_5+a_6y_6+a_7y_7) \ \xi \notin \operatorname{span}(y_{123},y_{147},y_{156},y_{257},y_{345},y_{367}).
\]
Indeed, the coefficient of $y_{246}$ in the expansion of the left-hand side is $-2a_2-2a_4-2a_6$. The relative Hodge--Riemann relations of Theorem~\ref{thmRHR} for $\operatorname{Gr}^{1,2}$ state that 
\[
(a_1y_1+a_2y_2+a_3y_3+a_4y_4+a_5y_5+a_6y_6+a_7y_7) \ \xi^2 
\]
is necessarily a negative multiple of $y_E$. 
A straightforward computation shows that the displayed element is, in fact, equal to $-8(a_2+a_4+a_6) y_E$. 
\end{example}

In general, in any decomposition of $M$ into indecomposable graded $B$-modules, the isomorphism classes and multiplicities of the summands are unique, but the summands need not be uniquely determined as submodules, as illustrated in Example~\ref{Fano1}. One implication of Theorem~\ref{thmHom} is that, when the graded $B$-module structure of $M$ is semismall,
$M$ has a unique isotypic decomposition into graded $B$-modules.

\begin{example}\label{Fano2}
Let $A$ be the graded M\"obius algebra of $F_7$ as above, 
and let $B$ be the subalgebra of $A$ generated by the elements $y_1$, $y_3+y_5$, and $y_2+y_4+y_6+y_7$. The $B$-module structure of $M$ is semismall because the sum of the three generators of $B$ is an element of $\mathscr{K}_A$. 
In this case, one can check that $M \coloneq A$ decomposes uniquely into indecomposable graded $B$-modules $N_1 \oplus N_2$, where
\begin{align*}
&N^3_1= \operatorname{span}(y_{E}), &N^3_2&=0,\\
&N^2_1=\operatorname{span}(y_{147},y_{246},y_{345},y_{123}+y_{156},y_{257}+y_{367}), &N^2_2&= \operatorname{span}(y_{123}-y_{156},y_{257}-y_{367}),\\
&N^1_1=\operatorname{span}(y_1,y_4,y_7,y_3+y_5,y_2+y_6), &N^1_2&=\operatorname{span}(y_3-y_5,y_2-y_6),\\
&N^0_1=\operatorname{span}(y_{\varnothing}), &N_2^0&=0.
\end{align*}
Note that the summands $N_1$ and $N_2$ are orthogonal with respect to $\mathscr{Q}$, as predicted by Theorems~\ref{thmHom} and ~\ref{thmPD}.
That $N_2$ is an indecomposable graded $B$-module follows from the property that 
\[
N_2^2=\operatorname{span}(y_1,y_3+y_5,y_2+y_4+y_6+y_7) \cdot \xi \ \ \text{for any nonzero $\xi \in N_2^1$.}
\]
That $N_2$ satisfies the Hodge--Riemann relations with respect to $\mathscr{K}_B$ is equivalent to the statement
\[
\text{$\begin{pmatrix}
a+3c& a-c\\
a-c&a+2b+c
\end{pmatrix}$ is positive definite for any $a,b,c>0$,}
\]
as can be checked directly.  In Example~\ref{Fano1}, each indecomposable Lefschetz module appearing in the decomposition of $M$ had a lowest-degree nonzero homogeneous component of dimension $1$.
In the current example, we see that indecomposable Lefschetz modules with lowest-degree components of dimension $2$  appear in the decomposition of $M$.
\end{example}

\begin{example}\label{VariousC}
Consider the graded real vector space $M=M^0 \oplus M^1 \oplus M^2$ with
\[
M^0=\mathbb{R}e^0_1 \oplus \mathbb{R} e^0_2, \quad
M^1=\mathbb{R}e^1_1 \oplus \mathbb{R} e^1_2, \quad
M^2=\mathbb{R}e^2_1 \oplus \mathbb{R} e^2_2.
\]
Let $\mathscr{Q}\colon M \times M \to \mathbb{R}$ be 
the symmetric bilinear form  satisfying 
\[
\mathscr{Q}(e^0_1,e^2_1)=\mathscr{Q}(e^1_1,e^1_1)=
\mathscr{Q}(e^0_2,e^2_2)=\mathscr{Q}(e^1_2,e^1_2)=1, \ \ \mathscr{Q}(e^2_1,e^0_2)=\mathscr{Q}(e^1_1,e^1_2)=\mathscr{Q}(e^0_1,e^2_2)=0,
\]
and the degree orthogonality relations 
$\mathscr{Q}(M^i,M^j)=0$ for $i+j\neq 2$. 
Let $A$ be the graded algebra generated by the commuting endomorphisms $y,z,z'$ of $M$ given by
\[
\begin{alignedat}{3}
y\cdot e^0_1 &= e^1_1, &\qquad z\cdot e^0_1 &= e^2_2, &\qquad z'\cdot e^0_1 &= e^2_1, \\
y\cdot e^1_1 &= e^2_1, & z\cdot e^1_1 &= 0, & z'\cdot e^1_1 &= 0,\\
y\cdot e^2_1 &= 0, & z\cdot e^2_1 &= 0, & z'\cdot e^2_1 &= 0, \\
y\cdot e^0_2 &= e^1_2, & z\cdot e^0_2 &= e^2_1, & z'\cdot e^0_2 &= -e^2_2, \\
y\cdot e^1_2 &= e^2_2, & z\cdot e^1_2 &= 0, & z'\cdot e^1_2 &= 0,\\
y\cdot e^2_2 &= 0, & z\cdot e^2_2 &= 0, & z'\cdot e^2_2 &= 0,
\end{alignedat}
\]
so that $A=A^0 \oplus A^1 \oplus A^2=\mathbb{R} 1 \oplus \mathbb{R}y \oplus (\mathbb{R} y^2 \oplus \mathbb{R} z \oplus \mathbb{R}z')$. Then $\mathscr{Q}$ is $A$-invariant, and $(M,\mathscr{Q})$ is a Lefschetz module of degree $2$ over $(A,\mathscr{K}_A)$, where $\mathscr{K}_A$ is the cone of positive multiples of $y$ in $A^1$.

Let $B$ be the graded subalgebra of $A$ generated by $y$, so that   the perverse filtration of $M$ over $B$ satisfies $P_0=P_1=0$ and $P_2=P_3=P_4=M$. 
Theorem~\ref{thmDecomposition} applies to any graded $B$-subalgebra $C$ of $A$. We note that, as $C$ varies, the conclusions of Theorem~\ref{thmDecomposition} are generally incomparable in strength.

\begin{enumerate}[(1)]\itemsep 5pt
\item When $C$ is generated by $y$, the graded $C$-module $M$ admits infinitely many decompositions of the form $N_1 \oplus N_2$, where $N_1 \simeq N_2$ as graded $C$-modules. The space of $C$-invariant symmetric bilinear forms on each $N_i$ satisfying the degree orthogonality relations is $1$-dimensional, as predicted by Theorem~\ref{thmDecomposition}. The space of $C$-invariant symmetric bilinear forms on $M$ satisfying the degree orthogonality relations is $3$-dimensional, and $\operatorname{Hom}_C(M,M)\simeq\operatorname{Mat}_{2\times 2}(\mathbb{R})$.
\item
When $C$ is generated by $y$ and $z$, the graded $C$-module $M$ has a unique decomposition $L_1 \oplus L_2$, where $L_1$ is generated by $e^0_1+e^0_2$ and $L_2$ is generated by $e^0_1-e^0_2$.
The space of $C$-invariant symmetric bilinear forms on each $L_i$ satisfying the degree orthogonality relations is $1$-dimensional, as predicted by Theorem~\ref{thmDecomposition}. 
The space of $C$-invariant symmetric bilinear forms on $M$ satisfying the degree orthogonality relations is $2$-dimensional, and $\operatorname{Hom}_C(M,M)\simeq \mathbb{R} \oplus \mathbb{R}$.
\item
When $C$ is generated by $y$ and $z'$, the graded $C$-module $M$ has a unique decomposition $L'_1 \oplus L'_2$, where $L'_1$ is generated by $e^0_1$ and $L'_2$ is generated by $e^0_2$.
The space of $C$-invariant symmetric bilinear forms on each $L'_i$ satisfying the degree orthogonality relations is $1$-dimensional, as predicted by Theorem~\ref{thmDecomposition}. 
The space of $C$-invariant symmetric bilinear forms on $M$ satisfying the degree orthogonality relations is $2$-dimensional, and $\operatorname{Hom}_C(M,M)\simeq \mathbb{R} \oplus \mathbb{R}$.
\item
When $C$ is generated by $y$ and $z$ and $z'$, the graded $C$-module $M$ is indecomposable. 
The space of $C$-invariant symmetric bilinear forms on $M$ satisfying the degree orthogonality relations is $1$-dimensional, as predicted by Theorem~\ref{thmDecomposition}, and  $\operatorname{Hom}_C(M,M)\simeq \mathbb{R}$.
\end{enumerate}
\end{example}

\begin{example}\label{ExampleHirzeburch}
According to \cite[Proposition 1.2.2]{KK87}, the set of elements of $A^1$ that satisfy the hard Lefschetz property and the Hodge--Riemann relations on $(M,\mathscr{Q})$ is a disjoint union of open convex cones $\mathscr{K}_1 \cup \cdots \cup \mathscr{K}_m$ in $A^1$,
referred to as the \emph{polarizing cones} of $(M,\mathscr{Q})$  in \cite{Cat08}.

For example, consider the graded algebra $A =A^0 \oplus A^1 \oplus A^2=\mathbb{R} y_\varnothing \oplus (\mathbb{R} y_1 \oplus \mathbb{R}y_2) \oplus \mathbb{R}y_E$, with the multiplicative structure given by
\[
y_\varnothing=1, \ \  y_1^2=y_2^2=0, \ \ \text{and} \ \ y_1y_2=y_E.
\]
Let  $\mathscr{Q}\colon A \times A \to \mathbb{R}$ be the unique $A$-invariant symmetric bilinear form satisfying  $\mathscr{Q}(1,y_E)=1$  and $\mathscr{Q}(A^i,A^j)=0$ for $i+j\neq 2$. 
Then $(A,\mathscr{Q})$ has two polarizing cones 
\[
\mathscr{K}_+=\{a_1y_1+a_2y_2 \mid a_1>0, a_2>0\} \ \ \text{and} \ \ \mathscr{K}_-=\{a_1y_1+a_2y_2 \mid a_1<0, a_2<0\}. 
\]
When $B$ is generated by $\ell \in \mathscr{K}_\pm$, the  $B$-module $A$ decomposes uniquely into  $L \oplus N$, where 
\[
L=B \cdot y_\varnothing \ \ \text{and} \ \ N=\ker(A^1 \xrightarrow{\ell} A^2).
\]
When $B$ is generated by $ y_1$, any decomposition of the graded $B$-module $A$ is of the form $L(t) \oplus N_1(t)$ for a real parameter $t$, where 
\[
L(t)=B\cdot y_\varnothing \ \ \text{and} \ \ N_1(t)=B \cdot (ty_1+y_2).
\]
When $B$ is generated by $ y_2$, any decomposition of  the graded $B$-module $A$ is of the form $L(t) \oplus N_2(t)$ for a real parameter $t$, where 
\[
L(t)=B\cdot y_\varnothing \ \ \text{and} \ \ N_2(t)=B \cdot (y_1+ty_2).
\]
Since $P_0=0$, we see from \eqref{eq:perverse-filtration} that $P_1=L(t)$ does not depend on $t$.

Note that the even cohomology of the $n$-th Hirzebruch surface $\mathbb{F}_n$ with real coefficients equipped with the Poincar\'e pairing is isomorphic to $(A,\mathscr{Q})$ over $\mathbb{R}$ in exactly two ways, up to scaling  $y_1$ and $y_2$. Under any such isomorphism, the ample cone of $\mathbb{F}_n$ maps to one of the polarizing cones, and the inclusion of the ample cone into the corresponding polarizing cone is necessarily proper if $n$ is positive. In this sense, the ample cone of $\mathbb{F}_n$ shrinks as $n$ increases.\footnote{This example also shows that \cite[Conjecture 2.1.7]{KK87} is not true as stated. For example, for the positive definite quadratic form $f(x,y)=6x^2-8xy+3y^2$, we have $f(y_1+y_2,y_1+2y_2)=0$ in $A$. One can also produce a counterexample to \cite[Conjecture 2.1.7]{KK87} based on \cite[Example 9.2]{RossToma}.} 
\end{example}

While the initial examples of Lefschetz modules came from algebraic geometry, there are Lefschetz modules which provably do not appear inside the cohomology of a compact K\"{a}hler manifold or the Chow ring of a smooth projective variety modulo numerical equivalence. 

\begin{example}\label{exReverseKT}
Consider the cubic polynomial in three variables 
\[
f(w_1, w_2, w_3) = 14w_1^3 + 6w_1^2w_2 +  24w_1^2 w_3 +  12w_1w_2w_3 + 6w_1w^2
_3 + 3w_2w_3^2,
\]
and let $A$ be the graded algebra over $\mathbb{R}$ cogenerated by $f$. In other words, $A$ is the unique quotient of $\mathbb{R}[x_1, x_2, x_3]$ equipped with an isomorphism $\deg \colon A^3 \simeq \mathbb{R}$ such that the bilinear forms $A^0 \times A^3 \to A^3 \simeq \mathbb{R}$ and $A^1 \times A^2 \to A^3 \simeq \mathbb{R}$ are nondegenerate and \[
\deg((w_1x_1 + w_2x_2 + w_3x_3)^3) = f(w_1, w_2, w_3).
\]
Using the fact that $f$ is a \emph{Lorentzian polynomial} in the sense of \cite{BHb}, one can check that $A$ is a Lefschetz module over $(A, \mathscr{K}_A)$, where $\mathscr{K}_A$ is the cone of positive linear combinations of $x_1,x_2,x_3$. However, $A$ cannot be a subquotient of the ring of real $(p,p)$-forms on a $3$-dimensional compact K\"{a}hler manifold nor of the Chow ring modulo numerical equivalence of a $3$-dimensional smooth projective variety over an algebraically closed field because $f$ does not satisfy the \emph{reverse Khovanskii--Teissier inequality} of \cite{LehmannXiao,JiangLi}. For more details, see \cite[Example 14]{HuhICM2022}. 
\end{example}

By Theorem~\ref{thmDecomposition}, the endomorphism ring of an indecomposable Lefschetz module is a finite-dimensional division algebra over $\mathbb{R}$. We give examples illustrating that both the field of complex numbers $\mathbb{C}$ and the quaternions $\mathbb{H}$ can occur.

\begin{example}\label{ex:endC}
Let $M^0$ be the $4$-dimensional real vector space with basis $e_1, e_2, e_3, e_4$, and let $M^1$ be the $4$-dimensional real vector space with basis $e_1^*, e_2^*, e_3^*, e_4^*$. Set $M \coloneq M^0 \oplus M^1$, and let $\mathscr{Q}\colon M \times M \to \mathbb{R}$ be the symmetric bilinear form defined by 
\[
\mathscr{Q}(e_i, e_j)=0, \ \ \mathscr{Q}(e_i^*, e_j^*)=0, \ \ \text{and} \ \ 
\mathscr{Q}(e_i, e_j^*) = \begin{cases} 1, & \text{if }i= j, \\ 0, & \text{otherwise}. \end{cases}
\]
A linear map $M^0 \to M^1$ is given by a $4 \times 4$ matrix, and $\mathscr{Q}$ is invariant under this map if and only if this matrix is symmetric. Set $A^0\coloneq\mathbb{R}$ and
\[
A^1 \coloneq \left \{ \begin{pmatrix}X & dY \\ dY^{t}& X \end{pmatrix} \;\middle|\;  X = \begin{pmatrix} a & c \\ c& b \end{pmatrix}, \ Y= \begin{pmatrix} 0 &  1 \\ -1& 0 \end{pmatrix}, \  a, b, c, d \in \mathbb{R} \right\}.
\]
The direct sum $A \coloneq A^0 \oplus A^1$ has the structure of a graded algebra over $\mathbb{R}$. Note that $M$ is an $A$-module and that $\mathscr{Q}$ is an $A$-invariant symmetric bilinear form on $M$. Setting $\mathscr{K}_A$ to be the intersection of $A^1$ with the cone of $4 \times 4$ positive definite matrices gives $(M, \mathscr{Q})$ the structure of a Lefschetz module of degree one over $(A, \mathscr{K}_A)$. 

The ring of graded $A$-module endomorphisms of $M$ is identified with the space of $4 \times 4$ matrices that commute with all matrices in $A^1$. It is easy to check that this space of matrices is spanned by the identity and the matrix 
\[
J \coloneq \begin{pmatrix} 0_2 & I_2 \\ -I_2 & 0_2  
\end{pmatrix},  \ \ \text{where} \ \ 0_2 = \begin{pmatrix} 0 & 0 \\ 0& 0 \end{pmatrix} \ \ \text{and} \ \ I_2= \begin{pmatrix} 1 & 0 \\ 0& 1 \end{pmatrix}.
\]
As the algebra spanned by the identity and $J$ is isomorphic to the complex numbers $\mathbb{C}$, we see that $M$ is an indecomposable $A$-module and $\operatorname{Hom}_A(M,M) \simeq \mathbb{C}$. 
In this case, the space of $A$-invariant bilinear forms on $M$ satisfying the degree orthogonality relations is $2$-dimensional, whereas the subspace of $A$-invariant symmetric bilinear forms satisfying the degree orthogonality relations is $1$-dimensional, as predicted by Theorem~\ref{thmPD}.
\end{example}

\begin{example}\label{ex:endH}
Let $M^0$ be the $8$-dimensional real vector space with basis $e_1, e_2, \dotsc, e_8$, and let $M^1$ be the $8$-dimensional real vector space with basis $e_1^*, e_2^*, \dotsc, e_8^*$. Set $M \coloneq M^0 \oplus M^1$, and let $\mathscr{Q}\colon M \times M \to \mathbb{R}$ be the symmetric bilinear form defined by 
\[
\mathscr{Q}(e_i, e_j)=0, \ \ \mathscr{Q}(e_i^*, e_j^*)=0, \ \ \text{and} \ \ \mathscr{Q}(e_i, e_j^*) = \begin{cases} 1, & \text{if }i= j, \\ 0, & \text{otherwise}. \end{cases}
\]
A linear map $M^0 \to M^1$ is given by an $8 \times 8$ matrix, and $\mathscr{Q}$ is invariant under this map if and only if this matrix is symmetric. Set $A^0 \coloneq \mathbb{R}$ and
\[
A^1 \coloneq \left \{ \begin{pmatrix} X & eY & dY & f Y\\
e Y^t & X & fY & dY^t \\ 
d Y^t & f Y^t & X & e Y\\ 
f Y^t & d Y & e Y^t & X
\end{pmatrix} \;\middle|\;  X= \begin{pmatrix}a&c\\ c&b\end{pmatrix}, \ Y= \begin{pmatrix}0&1\\ -1&0\end{pmatrix}, \ a, b, c, d, e, f \in \mathbb{R} \right\}.
\]
The direct sum $A \coloneq A^0 \oplus A^1$ has the structure of a graded algebra over $\mathbb{R}$. 
Let $\mathscr{K}_A$ be the intersection of $A^1$ with the cone of $8 \times 8$ positive definite matrices. Note that $(M, \mathscr{Q})$ is a Lefschetz module of degree one over $(A, \mathscr{K}_A)$. A direct computation shows that the space of $8 \times 8$ matrices that commute with $A^1$ is spanned by the identity, $J$, $K$, and $JK$, where
\[J \coloneq \begin{pmatrix} 0_2 & 0_2 & I_2 & 0_2 \\ 
  0_2 & 0_2 & 0_2 & I_2 \\ 
  - I_2 & 0_2 & 0_2 & 0_2\\ 
  0_2 & -I_2 & 0_2 & 0_2 
\end{pmatrix}, \ \ \text{and} \ \   K \coloneq \begin{pmatrix} 
  0_2 & 0_2 & 0_2 & I_2 \\
  0_2 & 0_2 & - I_2 & 0_2 \\ 
  0_2 & I_2 & 0_2 & 0_2 \\
  -I_2 & 0_2 & 0_2 & 0_2 
\end{pmatrix}.
\]
As the algebra spanned by the identity, $J$, $K$, and $JK$ is isomorphic to the quaternions $\mathbb{H}$,  
we see that $M$ is an indecomposable $A$-module and $\operatorname{Hom}_A(M,M) \simeq \mathbb{H}$. 
In this case, the space of $A$-invariant bilinear forms on $M$ satisfying the degree orthogonality relations is $4$-dimensional, whereas the subspace of $A$-invariant symmetric bilinear forms satisfying the degree orthogonality relations is $1$-dimensional, as predicted by Theorem~\ref{thmPD}.
\end{example}

\begin{remark}
Since an indecomposable Lefschetz module $M$ is a graded vector space over its endomorphism ring, Theorem~\ref{thmHom} implies that the endomorphism ring of such a module is necessarily isomorphic to $\mathbb{R}$ if some graded piece $M^k$ is an odd-dimensional vector space over $\mathbb{R}$. 
\end{remark}

The next examples show that the uniqueness of $\mathscr{Q}_\alpha$ in Theorems~\ref{thmPD} and ~\ref{thmDecomposition} may fail if one works with coefficients in $\mathbb{Q}$ instead of $\mathbb{R}$. 

\begin{example}\label{ex:Q-nonuniqueness}
Let $M^0$ be the $2$-dimensional rational vector space with basis $e_1, e_2$, and let $M^1$ be the $2$-dimensional rational vector space with basis $e_1^*, e_2^*$. Set $M \coloneq M^0 \oplus M^1$, and let $\mathscr{Q}\colon M \times M \to \mathbb{Q}$ be the symmetric bilinear form defined by 
\[
\mathscr{Q}(e_i, e_j)=0, \ \ \mathscr{Q}(e_i^*, e_j^*)=0, \ \ \text{and} \ \ 
\mathscr{Q}(e_i, e_j^*) = \begin{cases} 1, & \text{if }i= j, \\ 0, & \text{otherwise}. \end{cases}
\]
A linear map $M^0 \to M^1$ is given by a $2 \times 2$ matrix, and $\mathscr{Q}$ is invariant under this map if and only if this matrix is symmetric. Set $A^0\coloneq\mathbb{Q}$ and
\[
A^1 \coloneq \left \{ \begin{pmatrix}a & b \\ b& a + 2b \end{pmatrix} \;\middle|\; a, b \in \mathbb{Q}\right\}.
\]
The direct sum $A \coloneq A^0 \oplus A^1$ has the structure of a graded algebra over $\mathbb{Q}$. Note that $M$ is an $A$-module and that $\mathscr{Q}$ is an $A$-invariant symmetric bilinear form on $M$. Set $A_\mathbb{R} \coloneq A \otimes_{\mathbb{Q}} \mathbb{R}$, $M_\mathbb{R} \coloneq M \otimes_{\mathbb{Q}} \mathbb{R}$, $\mathscr{Q}_\mathbb{R} \coloneq \mathscr{Q} \otimes_{\mathbb{Q}} \mathbb{R}$, and let
$\mathscr{K}_A$ be the intersection of $A^1_\mathbb{R}$ with the cone of $2 \times 2$ positive definite matrices. Then $(M_\mathbb{R}, \mathscr{Q}_\mathbb{R})$ is a Lefschetz module of degree one over $(A_\mathbb{R}, \mathscr{K}_A)$. 

The ring of graded $A$-module endomorphisms of $M$ is identified with the space of $2 \times 2$ matrices that commute with all matrices in $A^1$. It is easy to check that this space of matrices is spanned by the identity and the matrix 
\[
D \coloneq \begin{pmatrix} -1 & 1 \\ 1 & 1  
\end{pmatrix}.
\]
Thus $\operatorname{Hom}_A(M,M) \simeq \mathbb{Q}[\sqrt{2}]$, and hence $M$ is an indecomposable graded $A$-module. The space of $A$-invariant symmetric bilinear forms $\mathscr{Q}_\alpha$ on $M$ satisfying the degree orthogonality relations is the $2$-dimensional space over $\mathbb{Q}$ spanned by
\[
(x_1, x_2) \longmapsto \mathscr{Q}(x_1, x_2)
\ \ \text{and} \ \ (x_1, x_2) \longmapsto \mathscr{Q}(Dx_1, x_2).
\] 
Some of these $\mathscr{Q}_\alpha$ satisfy the Hodge--Riemann relations up to a sign, whereas others do not. 
Note that $M_\mathbb{R}$ is a decomposable graded $A_\mathbb{R}$-module, and each of its indecomposable summands admits an essentially unique nonzero real-valued $A$-invariant symmetric bilinear form satisfying the degree orthogonality relations.
\end{example}

\begin{example}\label{ex:Q-nonuniqueness2}
Consider the graded $\mathbb{Q}$-algebra with trivial grading 
\[
A=A^0=\mathbb{Q}[x]/(x^2-x-1),
\]
and consider  the graded $A$-module with trivial grading $M=M^0=A^0$. 
Let $\mathscr{Q} \colon M \times M \to \mathbb{Q}$ be
the  $A$-invariant symmetric bilinear form  given by
\[
\mathscr{Q}(1,1)=1, \ \ \mathscr{Q}(x,x)=1, \ \ \text{and} \ \ \mathscr{Q}(1,x)=0.
\]
Since $\mathscr{Q}$ is positive definite, 
 $(M_\mathbb{R},\mathscr{Q}_\mathbb{R})$ is a Lefschetz module of degree $0$ over $(A_\mathbb{R},\mathscr{K}_A)$, where $\mathscr{K}_A$ is the zero cone in the zero vector space $A^1_\mathbb{R}$. 
 Since $A_\mathbb{R}$ preserves the perverse filtration of $M_\mathbb{R}$ over $\mathbb{R}$, which is trivial, Theorem~\ref{thmDecomposition} applies to the $A_\mathbb{R}$-module $M_\mathbb{R}$.

 The $A$-module $M$ is indecomposable, and the space of $A$-invariant symmetric bilinear forms $\mathscr{Q}_\alpha \colon M \times M \to \mathbb{Q}$ is $2$-dimensional  over $\mathbb{Q}$. Some of these $\mathscr{Q}_\alpha$ satisfy the Hodge--Riemann relations  up to a sign, whereas others do not.
On the other hand, we have a unique decomposition of $M_\mathbb{R}$ into two $1$-dimensional $A_\mathbb{R}$-submodules $L_1 \oplus L_2$ with $L_1 \nsimeq L_2$,
and  the space of $A_\mathbb{R}$-invariant symmetric bilinear forms $L_i \times L_i \to \mathbb{R}$ is $1$-dimensional  over $\mathbb{R}$ for each $i$.
\end{example}

One ingredient in the proof of the uniqueness of $\mathscr{Q}_\alpha$ in Theorem~\ref{thmDecomposition} is that, when there are two independent such forms on $M_\mathbb{R}$, there must be a nonzero degenerate $A_\mathbb{R}$-invariant symmetric bilinear form $M_\mathbb{R} \times M_\mathbb{R} \to \mathbb{R}$. In the settings of Examples~\ref{ex:Q-nonuniqueness} and ~\ref{ex:Q-nonuniqueness2}, it is easy to see that no such form is defined over $\mathbb{Q}$.

\section{The decomposition package in practice}\label{sec:applications}
In this section, we discuss some variations on the main theorems and give applications.

\subsection{Matroids}\label{sec:matroids}

In recent years, significant progress in matroid theory has been made by proving that certain combinatorially defined modules are Lefschetz; see, for example, \cite{AHK,ADH1,BHMPW2}. 
Let $\mathrm{M}$ be a simple matroid of rank $d$ on a finite ground set $E$. For each $0 \le k \le d$, let $\mathcal{L}^k(\mathrm{M})$ denote the set of rank $k$ flats of $\mathrm{M}$. 

\begin{theorem}\label{thm:topheavy}
For any nonnegative $k \le j \le d - k$, we have $|\mathcal{L}^k(\mathrm{M})| \le |\mathcal{L}^j(\mathrm{M})|$. 
\end{theorem}

Theorem~\ref{thm:topheavy} was conjectured by Dowling and Wilson in \cite{DW1,DW2}. It was proved for realizable matroids in \cite{HW} using the hard Lefschetz theorem for $\ell$-adic intersection cohomology of projective varieties over finite fields, and in full generality in \cite{BHMPW2} through the construction of the intersection cohomology of a matroid and the proof of the hard Lefschetz theorem for it. 
Its best-known special case is the de Bruijn--Erd\H{o}s theorem on point-line incidences in projective planes \cite{dBE}: 
\begin{quote}
\emph{Let $E$ be a set of points in a projective plane that is not contained in any line. Then there are at least $|E|$ distinct lines that intersect $E$ in at least two points.}
\end{quote}
We now show how the main results of \cite{BHMPW2}, and in particular Theorem~\ref{thm:topheavy}, can be deduced from the main results of this paper.

Let $\operatorname{H}(\mathrm{M})$ denote the \emph{graded M\"{o}bius algebra} of $\mathrm{M}$ over $\mathbb{R}$. As a graded vector space, 
\[
\operatorname{H}(\mathrm{M}) = \bigoplus_{k=0}^{d} \Bigg(\bigoplus_{F \in \mathcal{L}^k(\mathrm{M})} \mathbb{R}y_F\Bigg), \ \ \text{with $\operatorname{deg}(y_F)=\operatorname{rk}(F)$.}
\]
The multiplicative structure is given by the formula
\[
y_F  y_G = \begin{cases} y_{F \vee G} & \text{if }\operatorname{rk}(F) + \operatorname{rk}(G) = \operatorname{rk}(F \vee G), \\ 0 & \text{otherwise}. \end{cases}
\]
See Examples~\ref{Fano1} and \ref{Fano2} for the case of the Fano matroid. 
Note that $\operatorname{H}(\mathrm{M})$ is generated as an $\mathbb{R}$-algebra by the elements $y_i$ for each rank $1$ flat $i$.  
In \cite{BHMPW}, the authors construct a graded algebra $\operatorname{CH}(\mathrm{M})$, called the \emph{augmented Chow ring} of $\mathrm{M}$, which contains $\operatorname{H}(\mathrm{M})$ as a subalgebra. By \cite[Theorem 1.3]{BHMPW}, the augmented Chow ring is equipped with a nonempty open convex cone $\mathscr{K}$ in its degree one part and a linear map $\deg \colon \operatorname{CH}(\mathrm{M}) \to \mathbb{R}$ such that the symmetric bilinear form  $\mathscr{Q}(x_1, x_2) = \deg(x_1x_2)$ gives $(\operatorname{CH}(\mathrm{M}), \mathscr{Q})$ the structure of a Lefschetz module of degree $d$ over $\operatorname{CH}(\mathrm{M})$ with respect to  $\mathscr{K}$. Furthermore, each $y_i$ lies in the closure of $\mathscr{K}$. 

Choose a decomposition of $\operatorname{CH}(\mathrm{M})$ into indecomposable graded $\operatorname{H}(\mathrm{M})$-modules. Let $\operatorname{IH}(\mathrm{M})$ denote the unique summand which intersects the degree $0$ part of $\operatorname{CH}(\mathrm{M})$, so $\operatorname{IH}(\mathrm{M})$ contains $\operatorname{H}(\mathrm{M})$. By the Krull--Schmidt theorem,  
the graded  $\operatorname{H}(\mathrm{M})$-module $\operatorname{IH}(\mathrm{M})$ is well-defined up to isomorphism. 
This object is, by definition, the \emph{intersection cohomology} of $\mathrm{M}$ introduced in \cite{BHMPW2}.

\begin{proof}[Proof of Theorem~\ref{thm:topheavy}]
The interior of $\overline{\mathscr{K}} \cap \operatorname{H}^1(\mathrm{M})$ contains the set of all positive linear combinations of the $y_i$.
Thus, by Theorem~\ref{thmHL}, multiplication by $(\sum y_i)^{d - 2k}$ is an isomorphism from $\operatorname{IH}^k(\mathrm{M})$ to $\operatorname{IH}^{d-k}(\mathrm{M})$. In particular, multiplication by $(\sum y_i)^{d - 2k}$ induces an injection from $\operatorname{H}^k(\mathrm{M})$ to $\operatorname{H}^{d-k}(\mathrm{M})$, as multiplication by $\sum y_i$ preserves the submodule $\operatorname{H}(\mathrm{M})$ of $\operatorname{IH}(\mathrm{M})$. This injection factors through $\operatorname{H}^j(\mathrm{M})$ for each $k \le j \le d-k$. Since the dimension of $\operatorname{H}^k(\mathrm{M})$ is the cardinality of $\mathcal{L}^k(\mathrm{M})$, this implies the desired inequality. 
\end{proof}

Applied to a decomposition of the $\operatorname{H}(\mathrm{M})$-module $\operatorname{CH}(\mathrm{M})$, the decomposition package developed in this paper recovers Poincar\'e duality, the hard Lefschetz theorem, and the Hodge--Riemann relations for $\operatorname{IH}(\mathrm{M})$ established in \cite[Theorem 1.6]{BHMPW2}. It also yields the relative hard Lefschetz theorem and the relative Hodge--Riemann relations for $\operatorname{CH}(\mathrm{M})$ over $\operatorname{H}(\mathrm{M})$, which are new.

The decomposition package also applies to several other constructions associated with matroids. For example, in \cite{BHMPW}, the authors study the decomposition of the augmented Chow ring of a matroid $\mathrm{M}$ as a module over the corresponding ring of a deletion $\mathrm{M}\setminus S$ for $S \subseteq E$. Since the augmented Chow ring of the deletion is a subalgebra generated by elements of $\overline{\mathscr{K}}$, Theorems~\ref{thmRHL} and~\ref{thmRHR} apply to this decomposition, yielding the relative hard Lefschetz theorem and the relative Hodge--Riemann relations in this setting.

\subsection{Polytopes}\label{ssec:polytopes}

We now give an application of the main theorems to invariants of subdivisions of polytopes. In \cite{subdivisions}, Stanley introduced the \emph{local $h$-vector} of a triangulation of a simplex, and he used the relative hard Lefschetz theorem of \cite{BBD} to prove that local $h$-vectors of regular triangulations of simplices are unimodal. In \cite[page 134]{StanleyCCA}, Stanley asked for a proof of this result which does not use results from algebraic geometry. We explain how Stanley's result can be deduced from the main theorems of this paper and the work of McMullen on polytope algebras \cite{McMullen}. An alternative argument was given by Karu \cite{KaruRHL}. 

A \emph{geometric triangulation} of a simplex $\Delta$ is a simplicial complex $\widetilde \Delta$ together with an identification of the geometric realization of $\widetilde \Delta$ with the geometric realization of $\Delta$, where $\Delta$ and each simplex of $\widetilde \Delta$ are realized as Euclidean simplices. The triangulation is \emph{regular} if it arises from  projecting the lower faces of a convex polytope  obtained by lifting some of the points of $\Delta$ and taking the convex hull. 

Given a geometric triangulation of a $(d-1)$-dimensional simplex $\Delta$, we obtain a proper toric morphism to the affine space $\mathbb{A}^d$ as follows: View the geometric realization of $\Delta$ as the standard simplex in $\mathbb{R}^d$. 
The cone over the triangulation $\widetilde \Delta$ is then a fan $\Sigma$ supported on the positive orthant. By perturbing the vertices, we may assume that the fan is rational; if the initial triangulation was regular, then the perturbed triangulation will also be regular. The toric variety $X$ of   $\Sigma$ admits a proper toric morphism $f \colon X \to \mathbb{A}^d$, which  is projective  if and only if the triangulation is regular. 

Instead of giving the combinatorial definition of the local $h$-vector of a triangulation of a simplex \cite[Definition 2.1]{subdivisions}, we describe Stanley's toric interpretation of the local $h$-vector \cite[Theorem 5.2]{subdivisions}. 
The constant sheaf $\underline{\mathbb{R}}_{X}$ coincides with the intersection complex up to a shift because $X$ is a simplicial toric variety. Since $f$ is toric, when we apply the decomposition theorem \cite[Theorem 6.2.5]{BBD} to $f_* \underline{\mathbb{R}}_{X}$, all of the summands are supported on torus-orbit closures in $\mathbb{A}^d$. The local $h$-vector $(\ell_0, \dotsc, \ell_d)$ of the triangulation  counts the number of summands supported at the origin $o$ of $\mathbb{A}^d$ according to their shifts. More precisely, in any decomposition of $f_* \underline{\mathbb{R}}_{X}$,
 there are $\ell_i$ summands isomorphic to $\underline{\mathbb{R}}_o[-2i]$. In \cite[Theorem 5.2]{subdivisions}, Stanley deduced the following result from the relative hard Lefschetz theorem of \cite[Theorem 6.2.10]{BBD}. 

\begin{theorem}\label{thm:localhunimodal}
For any regular triangulation of a $(d-1)$-dimensional simplex, we have $\ell_i = \ell_{d-i}$ and $\ell_0 \le \ell_1 \le \dotsb  \le \ell_{\lfloor d/2 \rfloor}$. 
\end{theorem}

We explain how Theorem~\ref{thm:localhunimodal} follows from  Corollary~\ref{cor:multiplicities} and McMullen's work on polytope algebras \cite{McMullen}. 
Let $\overline{\Sigma}$ be the fan obtained by adding the ray generated by $(-1, \dotsc, -1)$ to the fan $\Sigma$ and, for each cone $\sigma$ in ${\Sigma}$ which is contained in the boundary of the positive orthant, adding the convex hull of $\sigma$ with this ray. The fan $\overline{\Sigma}$ is a subdivision of the fan of $\mathbb{P}^d$, so there is an induced map of the corresponding toric varieties $\overline{f}:\overline{X} \to \mathbb{P}^d$, which is projective if and only if the triangulation is regular. 
In particular, in this case, $\overline{X}$ is projective. 

The degree one part of the Chow ring $A$ of $\overline{X}$ is identified with the space of piecewise linear functions on $\overline{\Sigma}$ modulo the space of globally linear functions, and it is equipped with an open convex cone $\mathscr{K}$ given by the classes of strictly convex piecewise linear functions. The Chow ring is also equipped with a linear map $\deg \colon A^d \to \mathbb{R}$ such that the symmetric bilinear form $\mathscr{Q}(x_1,x_2) = \deg(x_1x_2)$ gives $(A, \mathscr{Q})$ the structure of a Lefschetz module of degree $d$ over $A$ with respect to $\mathscr{K}$. A combinatorial proof of this fact was given in \cite{McMullen}. 

\begin{proof}[Proof of Theorem~\ref{thm:localhunimodal}]
Let $\ell$ be the class of a nonzero convex piecewise linear function which is pulled back from the fan of $\mathbb{P}^d$. Since $\ell$ is convex, it lies in the closure of $\mathscr{K}$. Consider the graded vector space
\[
V
=
\operatorname{ann}_{A}(\ell)
\big/
\bigl(
\ell A
\cap
\operatorname{ann}_{A}(\ell)
\bigr).
\]
In \cite[Theorem 4.6]{subdivisions}, Stanley constructed a graded vector space whose Hilbert function is given by the local $h$-vector. It is straightforward to identify $V$ with this graded vector space; see \cite[Lemma 2.5]{LSlocalh}. 
Since the Hilbert function of $V$ counts the number of $1$-dimensional summands  in a decomposition of $A$ as a graded $\ell$-module,  Theorem~\ref{thm:localhunimodal} follows from Corollary~\ref{cor:multiplicities}. 
\end{proof}

In \cite{BBFK,BrLu}, an analogue of intersection cohomology for non-rational polytopes was constructed, and it was proved to be a Lefschetz module in \cite{Karu}. In \cite{KaruRHL}, a version of the relative hard Lefschetz theorem in this setting was proved. This result implies the unimodality of local $h$-vectors of regular polyhedral subdivisions \cite{KatzStapledon}. 
Our main theorems can be used to formally deduce the main result of \cite{KaruRHL} from \cite{Karu}. 

\subsection{Algebraic cycles}\label{sec:algebraiccycles}

Assuming Grothendieck's standard conjectures \cite{Grothendieck}, the ring of algebraic cycles modulo numerical equivalence on a smooth projective variety is a Lefschetz module; see Proposition~\ref{prop:standardlefschetz}. This statement is known unconditionally in certain cases, including smooth projective varieties over $\mathbb{C}$ for which the Hodge conjecture holds, as well as several classes of varieties over fields of arbitrary characteristic \cite{ItoInventiones}. Whenever it holds for a smooth projective variety $X$, our main theorems can be applied to any morphism from $X$ to a projective variety $Y$ by taking $B$ to be the subalgebra generated by the pullbacks of ample divisor classes on $Y$. In some cases, this yields instances of the standard conjectures for $Y$. This approach is closely related to the work of Corti and Hanamura \cite{Corti,CortiII}, who develop a version of the decomposition theorem for Chow groups under several additional conjectural assumptions.

Let $X$ be a smooth projective variety of dimension $d$ over an
algebraically closed field of arbitrary characteristic.\footnote{Throughout this paper, a variety over $k$ is by definition a reduced and irreducible scheme of finite type over $k$.} Let
$A_{\operatorname{num}}(X)$ denote its Chow ring of cycles modulo numerical
equivalence with real coefficients, equipped with the degree isomorphism
\[
\deg\colon A_{\operatorname{num}}^{d}(X)\xrightarrow{\sim}\mathbb{R}.
\]
We consider the following versions of standard conjecture A and the standard
conjecture of Hodge type \cite{Grothendieck}.

\begin{conjecture}\label{ConjectureA}
Let $\eta$ be the class of an ample divisor.
For any nonnegative $k\leq d/2$, multiplication by $\eta^{d-2k}$
induces an isomorphism
\[
\eta^{d-2k}\colon
A_{\operatorname{num}}^k(X)
\xrightarrow{\sim}
A_{\operatorname{num}}^{d-k}(X).
\]
\end{conjecture}

\begin{conjecture}\label{ConjectureH}
Let $\eta$ be the class of an ample divisor.
For any nonnegative $k\leq d/2$, the symmetric bilinear form
\[
A_{\operatorname{num}}^k(X)\times A_{\operatorname{num}}^k(X)
\longrightarrow \mathbb{R},
\qquad
(x_1,x_2)
\longmapsto
(-1)^k\deg\bigl(\eta^{d-2k}x_1x_2\bigr),
\]
is positive definite on the primitive subspace $\ker(\eta^{d-2k+1})$.
\end{conjecture}

For a smooth projective variety $X$, let $\mathscr{K}_X$ be the ample cone in $A^1_{\operatorname{num}}(X)$ given by positive real linear combinations of classes of ample divisors  \cite[Definition 1.3.11]{Lazarsfeld}. 
Let $\mathscr{Q}_X$ be the symmetric bilinear form on $A_{\operatorname{num}}(X)$ given by $(x_1, x_2) \mapsto \deg(x_1 x_2)$. 

\begin{proposition}\label{prop:standardlefschetz}
If Conjectures~\ref{ConjectureA} and ~\ref{ConjectureH} hold for all smooth projective varieties, then, for any $d$-dimensional smooth projective variety $X$,  
$(A_{\operatorname{num}}(X), \mathscr{Q}_X)$ is a Lefschetz module of degree $d$ over $(A_{\operatorname{num}}(X), \mathscr{K}_X)$. 
\end{proposition}

\begin{proof}
We argue by induction on the dimension $d \ge 2$ of $X$. The cases $d \le 1$ are immediate.

Let $\eta\in\mathscr{K}_X$, and suppose that, for some integer $k<d/2$, a
nonzero class $g\in A_{\operatorname{num}}^k(X)$ lies in the kernel of
multiplication by $\eta^{d-2k}$. Write
\[
\eta=\sum_i a_i[D_i],
\]
where each $a_i$ is a positive real number and each $D_i$ is an ample
divisor. By Bertini's theorem, we may choose the divisors $D_i$ to be smooth
and irreducible. We then have
\[
0
=
\mathscr{Q}_X\bigl(g,\eta^{d-2k}g\bigr)
=
\sum_i a_i\,
\mathscr{Q}_X\bigl(g,[D_i]\eta^{d-2k-1}g\bigr).
\]
Let $\iota_i^*$ be the restriction map $A_{\operatorname{num}}(X) \to A_{\operatorname{num}}(D_i)$. 
By the projection formula, we have 
\[
\mathscr{Q}_X(g, [D_i]\eta^{d-2k-1} g) = \mathscr{Q}_{D_i}(\iota_i^*(g), \iota_i^*(\eta^{d-2k-1} g)).
\] 
We claim that $\iota_i^*(g)$ is nonzero for each $i$. 
Indeed, Conjecture~\ref{ConjectureA} implies that $[D_i]^{d - 2k} g$ is nonzero. Hence there is a class $h_i \in A^{k}_{\operatorname{num}}(X)$ such that 
\[
\mathscr{Q}_X(h_i, [D_i]^{d - 2k}g) = \mathscr{Q}_{D_i}(\iota_i^*(h_i), \iota_i^*([D_i]^{d - 2k - 1})\iota_i^*(g)) \neq 0.
\]
Since $\iota^*_i(\eta^{d - 2k} g) = 0$ and $A_{\operatorname{num}}(D_i)$ is a Lefschetz module by induction, the Hodge--Riemann relations 
imply that 
\[
(-1)^k\mathscr{Q}_X(g, [D_i]\eta^{d-2k-1} g) > 0.
\]
Applying this inequality for every $i$ contradicts  $\mathscr{Q}_X(g, \eta^{d-2k}g) = 0$. 
This verifies the hard Lefschetz property for any $\eta \in \mathscr{K}_X$. 
It follows that the signature of the Hodge--Riemann form on $A^k_{\operatorname{num}}(X)$  does not change as we vary $\eta$ within $\mathscr{K}_X$. Recall that the Hodge--Riemann relations 
for $\eta$ are equivalent to the signature of this form being $\sum_{i=0}^{k} (-1)^{i} (\dim A^i_{\operatorname{num}}(X) - \dim A^{i-1}_{\operatorname{num}}(X))$ for each $k$. 
By Conjecture~\ref{ConjectureH}, this holds when $\eta$ is the class of an ample divisor; hence the Hodge--Riemann relations hold for every $\eta \in \mathscr{K}_X$.
\end{proof}

That $A_{\operatorname{num}}(X)$ is a Lefschetz module is sufficient for many applications of the standard conjectures. For example, if $X$ is a smooth projective variety over the algebraic closure of a finite field, then $A_{\operatorname{num}}(X \times X)$ being a Lefschetz module implies the Riemann hypothesis for $X$ \cite{Kleiman}. 

Let $A$ be a graded $\mathbb{R}$-algebra equipped with an isomorphism $\deg \colon A^d \to \mathbb{R}$ and a nonempty open convex cone $\mathscr{K}_A$ in $A^1$. Let $\mathscr{Q}_A$ be the symmetric bilinear form defined by $(x_1, x_2) \mapsto \deg(x_1x_2)$. 
Suppose that $(A, \mathscr{Q}_A)$ is a Lefschetz module of degree $d$ over $(A, \mathscr{K}_A)$.

\begin{lemma}\label{lem:degpositive}
Let $y_1,\dotsc,y_e$ be elements of $\overline{\mathscr{K}}_A$. If the product 
$y_1\cdots y_e$ is nonzero, then
\[
\deg\bigl(\eta^{d-e}y_1\cdots y_e\bigr)>0 \ \ \text{for every $\eta \in \mathscr{K}_A$.}
\]
\end{lemma}

\begin{proof}
We iteratively apply Lemma~\ref{lem:descent} to conclude that  $A/\operatorname{ann}(y_1 \dotsb y_e)$, equipped with the bilinear form $\mathscr{Q}$ described in Section~\ref{ssec:descent}, is a nonzero Lefschetz module of degree $d-e$ over $(A, \mathscr{K}_A)$. 
The Hodge--Riemann relations 
for $k=0$ imply the required positivity. 
\end{proof}

Let $B$ be a graded subalgebra  of $A$ for which $B^e$ is one-dimensional and $B^s$ is zero-dimensional  for $s > e$. Choose an isomorphism $\deg_B \colon B^e \to \mathbb{R}$, and let $\mathscr{Q}_B$ be the symmetric bilinear form given by $\mathscr{Q}_B(x_1, x_2) = \deg_B(x_1x_2)$. 
We assume that $\mathscr{Q}_B$ is nondegenerate. This condition is independent of the choice of isomorphism $\deg_B$.

\begin{proposition}\label{prop:deg1}
Suppose that $B$ is generated as an $\mathbb{R}$-algebra by elements of
$\overline{\mathscr{K}}_A$, and let $\mathscr{K}_B$ be the interior of
$\overline{\mathscr{K}}_A\cap B^1$ in $B^1$. Then there is a unique sign
$\epsilon=\pm 1$ such that $(B,\epsilon\mathscr{Q}_B)$ is a Lefschetz
module of degree $e$ over $(B,\mathscr{K}_B)$.
\end{proposition}

\begin{proof}
Since $\mathscr{K}_B$ is open in $B^1$, we may choose elements $y_1, \dotsc, y_s$ from $\mathscr{K}_B$ that generate $B$ as an $\mathbb{R}$-algebra. Let $\eta$ be a class in $\mathscr{K}_A$, and set $\ell = y_1 + \dotsb + y_s$. Then
\[
\deg(\eta^{d-e} \ell^{e}) = \sum_{(i_1, \dotsc, i_e)} \deg(\eta^{d-e} y_{i_1} \dotsb y_{i_e}).
\]
By Lemma~\ref{lem:degpositive}, each term in this sum either vanishes 
or is positive. Since $B$ is generated by $y_1, \dotsc, y_s$ and $B^e \neq 0$, at least one term is positive. We deduce that $\ell^e \not \in \operatorname{ann}(\eta^{d-e})$. Since $B^e$ is
one-dimensional, it follows that
\[
B^e\cap\operatorname{ann}(\eta^{d-e})=0.
\]

The image of $B$ under the quotient map $A \to A/\operatorname{ann}(\eta^{d-e})$ is isomorphic to $B$. Indeed, the nondegeneracy of $\mathscr{Q}_B$ implies that every nonzero ideal in $B$ intersects $B^e$. As $B^e \cap \operatorname{ann}(\eta^{d-e}) = 0$, this implies that $B \cap \operatorname{ann}(\eta^{d-e}) = 0$. By Lemma~\ref{lem:descent}, $A/\operatorname{ann}(\eta^{d-e})$ is a Lefschetz module over $(A, \mathscr{K}_A)$. We may replace $A$ by $A/\operatorname{ann}(\eta^{d-e})$, and all the hypotheses are still satisfied. We may therefore assume that $d = e$. 

Up to a nonzero constant, the restriction of $\mathscr{Q}_A$ to $B$ is $\mathscr{Q}_B$, and so the restriction of $\mathscr{Q}_A$ to $B$ is nondegenerate. This implies that $B$ is a $B$-module summand of $A$, and it is  indecomposable as a graded $B$-module. Theorems~\ref{thmHL} and ~\ref{thmHR}  then show that $(B, \mathscr{Q}_A|_{B})$ is a Lefschetz module over $(B, \mathscr{K}_B)$. We set $\epsilon$ to be the sign of the constant relating $\mathscr{Q}_A|_{B}$ to $\mathscr{Q}_B$. 
\end{proof}

Let $X \to Y$ be a surjective map of smooth projective varieties, and suppose that $(A_{\operatorname{num}}(X), \mathscr{Q}_X)$ is a Lefschetz module of degree $\dim X$ over $(A_{\operatorname{num}}(X), \mathscr{K}_X)$. 

\begin{theorem}\label{thm:standarddeg1}
If $A_{\operatorname{num}}(Y)$ is generated by its degree one part as an $\mathbb{R}$-algebra, then $(A_{\operatorname{num}}(Y), \mathscr{Q}_Y)$ is a Lefschetz module of degree $\dim Y$ over $(A_{\operatorname{num}}(Y), \mathscr{K}_Y)$.
\end{theorem}

\begin{proof}
The hypotheses of Proposition~\ref{prop:deg1} hold, so it remains to check that the constant $\epsilon$ appearing there is $1$. This follows from the fact that $\mathscr{Q}_Y(1, \ell^{\dim Y}) > 0$ for  an ample class $\ell$ on $Y$. 
\end{proof}

\begin{remark}\label{Remark:LStapledon}
Theorem~\ref{thm:standarddeg1} is used in \cite{LPartida} to prove cases of Conjectures~\ref{ConjectureA} and ~\ref{ConjectureH} when $X$ is a projective bundle over a toric variety, and also to prove that certain algebras associated to matroids which resemble the cohomology rings of projective bundles are Lefschetz modules over themselves. This latter result is used in \cite{LStapledon} to prove an analogue of Fulton--Lazarsfeld positivity in \cite{FultonLaz} for matroids, which was conjectured implicitly in \cite[page 890]{Speyerg} and explicitly in \cite[Conjecture 9.13]{BF22}. 
\end{remark}

\subsection{Complex projective varieties}\label{sec:hodge-decomposition-theorem}

We explain how the results of this paper can be used to give alternative proofs of several celebrated results about the intersection cohomology of algebraic varieties and the Hodge-theoretic nature of the decomposition theorem. The main theorems give a purely algebraic version of the decomposition theorem. The results of \cite{BBD} can be used, in the setting of a map $X \to Y$ of complex projective varieties, to relate the algebraic decomposition to the usual decomposition theorem for perverse sheaves; see Proposition~\ref{prop:geometricperverse}. Using only the main theorems and results which were available when \cite{BBD} was written, this allows us to show that the geometrically defined perverse filtration on the cohomology of $X$ is by Hodge substructures and that the relative Hodge--Riemann relations hold. This in turn implies that the intersection cohomology of a projective variety carries a polarized pure Hodge structure. 
This result attracted considerable attention in the 1970s and 1980s \cite{Z1,Cheeger, Z2,HP85,CKS87,KK87}, and was finally proved using M. Saito's theory of mixed Hodge modules~\cite{Saito88,Saito90}.

\subsubsection{Lefschetz modules with Hodge structures}\label{rmkComplexCoefficients}  
 The results of this paper can be extended to modules with complex coefficients. Let $A$ and $\mathscr{K}_A$ be as in Section~\ref{sec:introduction}, and set $A_{\mathbb{C}}=A\otimes_{\mathbb{R}}\mathbb{C}$.  Let $\mathscr{M} = \bigoplus_{k\geq 0} \mathscr{M}^k$ be a finite-dimensional graded $A_{\mathbb{C}}$-module endowed with a Hermitian form $\mathscr{Q}\colon \mathscr{M}\times \mathscr{M} \to \mathbb{C}$ that is \emph{$A_\mathbb{C}$-invariant}:
\[
\mathscr{Q}(ax_1,x_2)=\mathscr{Q}(x_1,\bar ax_2) \ \ \text{for all $a \in A_{\mathbb{C}}$ and all $x_1,x_2 \in \mathscr{M}$.}
\]
We say that $(\mathscr{M},\mathscr{Q})$ is a \emph{complex Lefschetz module} of degree $d$ over $(A,\mathscr{K}_A)$ if it satisfies the K\"ahler package for $(\mathscr{Q},\mathscr{K}_A)$:
\begin{enumerate}\itemsep 5pt
\item[(PD)] The induced sesquilinear pairing between the graded pieces 
\[
\mathscr{Q} \colon \mathscr{M}^i \times \mathscr{M}^j \to \mathbb{C}
\]
is nondegenerate if $i+j=d$ and it is zero if $i+j \neq d$ (\emph{Poincar\'e duality}).
\item[(HL)] For each nonnegative $k \le \frac{d}{2}$ and $\eta \in \mathscr{K}_A$, the linear map
\[
\mathscr{M}^k \longrightarrow \mathscr{M}^{d-k}, \qquad x \longmapsto \eta^{d-2k}  x
\]
is an isomorphism of vector spaces (\emph{hard Lefschetz property}).
\item[(HR)] For each nonnegative $k \le \frac{d}{2}$ and  $\eta \in \mathscr{K}_A$, the Hermitian form
\[
\mathscr{M}^k \times \mathscr{M}^k \longrightarrow \mathbb{C}, \qquad (x_1,x_2) \longmapsto (-1)^k \mathscr{Q}\Big(x_1, \eta^{d-2k}  x_2\Big)
\]
is positive definite on the kernel of the linear map
\[
\mathscr{M}^k \longrightarrow \mathscr{M}^{d-k+1}, \qquad x \longmapsto \eta^{d-2k+1} x
\]
(\emph{Hodge--Riemann relations}).
\end{enumerate}
If $M$ is a Lefschetz module of degree $d$ over $(A,\mathscr{K}_A)$, then $M_\mathbb{C} = M \otimes_\mathbb{R} \mathbb{C}$ is a complex Lefschetz module of degree $d$ over $(A,\mathscr{K}_A)$ in the obvious way.

The proofs of the main theorems can be easily adapted to the setting of complex Lefschetz modules. Let $B$ and $\mathscr{K}_B$ be as in Section~\ref{sec:introduction}, and set $B_{\mathbb{C}}=B\otimes_{\mathbb{R}}\mathbb{C}$. We choose a decomposition 
\[
\mathscr{M} \simeq \bigoplus_\alpha \bigoplus_{k} \mathscr{N}_\alpha[-k]^{\oplus m(\alpha,k)},
\]
where $\mathscr{N}_\alpha=\mathscr{N}^0_\alpha \oplus \cdots \oplus \mathscr{N}^{d(\alpha)}_\alpha$ are indecomposable graded $B_\mathbb{C}$-modules satisfying 
\[
\mathscr{N}^0_\alpha \neq 0, \ \ \mathscr{N}^{d(\alpha)}_\alpha \neq 0, \ \ \text{and} \ \  \text{$\mathscr{N}_\alpha \nsimeq \mathscr{N}_\beta$ as graded $B_\mathbb{C}$-modules for all $\alpha \neq \beta$}.  
\]
Then Theorems~\ref{thmHom}, ~\ref{thmPD}, ~\ref{thmHL},  and~\ref{thmHR} hold for complex Lefschetz modules in the following sense.

\begin{theorem}\label{thmHomcomplex}
For each $\alpha$, we have an isomorphism of $\mathbb{C}$-algebras
\[
\operatorname{Hom}_{B_\mathbb{C}}(\mathscr{N}_\alpha,\mathscr{N}_\alpha) \simeq \mathbb{C}.
\]
If $\alpha=\beta$ and $k$ is strictly positive, or if $\alpha\neq\beta$ and
$
d(\alpha)\le d(\beta)+2k$, then 
\[
\operatorname{Hom}_{B_\mathbb{C}}(\mathscr{N}_\alpha,\mathscr{N}_\beta[-k])=0.
\]
Up to a nonzero complex multiple,  there is a unique nonzero $B_\mathbb{C}$-invariant sesquilinear form $\mathscr{Q}_\alpha$ on $\mathscr{N}_\alpha$ satisfying the degree orthogonality relation
\[
\mathscr{Q}_\alpha(\mathscr{N}_\alpha^i,\mathscr{N}_\alpha^j)=0 \ \ \text{unless $i+j=d(\alpha)$.} 
\]
Then there is a unique $\epsilon_{\alpha}$ of norm $1$ such that $\epsilon_\alpha \mathscr{Q}_\alpha$ gives $\mathscr{N}_\alpha$ the structure of a complex Lefschetz module of degree $d(\alpha)$ over $(B, \mathscr{K}_B)$.
\end{theorem}

The perverse filtration $0 \subseteq \mathscr{P}_0 \subseteq \mathscr{P}_1 \subseteq \cdots \subseteq \mathscr{P}_{2d}$ 
of $\mathscr{M}$ is defined as in (\ref{eq:perverse-filtration-annihilators}) by
\[
\mathscr{P}_j=\bigoplus_{d(\alpha)+2k \le j} \mathscr{N}_\alpha[-k]^{\oplus m(\alpha, k)}.
\]
It follows from Theorem~\ref{thmHomcomplex} that $\mathscr{P}_j$ does not depend on the chosen decomposition of $\mathscr{M}$, and the associated graded $\mathscr{G}\mathrm{r}$  is a graded $B_{\mathbb{C}}$-module and a graded $A_\mathbb{C}$-module  with \[
B_{\mathbb{C}}^k \, \mathscr{G}\mathrm{r}^{i,j} \subseteq \mathscr{G}\mathrm{r}^{i+k,j} \quad \text{and} \quad A_{\mathbb{C}}^k * \mathscr{G}\mathrm{r}^{i,j} \subseteq \mathscr{G}\mathrm{r}^{i+k,j+2k}.
\]
The Hermitian form $\mathscr{Q}$ on $\mathscr{M}$ induces a Hermitian form
\[
\underline{\mathscr{Q}} \colon \mathscr{G}\mathrm{r}^{i,j} \times \mathscr{G}\mathrm{r}^{d-i,2d-j} \longrightarrow \mathbb{C},
\]
and this form is both $A_{\mathbb{C}}$-invariant and $B_{\mathbb{C}}$-invariant:
\[
\underline{\mathscr{Q}}(a * x_1,x_2)=
\underline{\mathscr{Q}}(x_1,\bar{a} * x_2)
\ \  \text{and} \ \ 
\underline{\mathscr{Q}}(bx_1,x_2)=
\underline{\mathscr{Q}}(x_1,\bar{b}x_2) \ \  \text{for all $a\in A_{\mathbb{C}}$ and $b\in B_{\mathbb{C}}$}.
\]
With these modifications, the complex analogues of Theorems~\ref{thmRHL} and~\ref{thmRHR}, as well as Corollaries~\ref{cor:multiplicities}, ~\ref{cor:kernellefschetz}, 
 and~\ref{cor:kernellefschetz2} hold.

\begin{example}
Let $X$ be a connected compact K\"{a}hler manifold of complex dimension $d$. 
For each nonnegative integer $r$, the $r$-th cohomology with complex coefficients of $X$ is equipped with the Hodge decomposition 
\[
\operatorname{H}^r(X, \mathbb{C}) = \bigoplus_{p+q=r} \operatorname{H}^{p,q}(X).
\]
Let $A = \bigoplus_{k} \operatorname{H}^{2k}(X, \mathbb{R}) \cap \operatorname{H}^{k,k}(X)$ so that $A_{\mathbb{C}} = \bigoplus_{k} \operatorname{H}^{k,k}(X)$. 
The Poincar\'e pairing $\mathscr{R}$ gives $A$ the structure of a Lefschetz module over $(A,\mathscr{K}_A)$, where $\mathscr{K}_A$ is the K\"{a}hler cone in $\operatorname{H}^{1,1}(X, \mathbb{R})$. 
For each nonnegative integer $n$, let $\mathscr{M}_n$ be the graded $A_{\mathbb{C}}$-module defined by 
\[
\mathscr{M}_n = \bigoplus_{k} \mathscr{M}_n^k = \bigoplus_{k} \operatorname{H}^{n+k,k}(X). 
\]
 We endow $\mathscr{M}_n$ with the Hermitian form 
\[
\mathscr{Q}_n \colon \mathscr{M}_n \times \mathscr{M}_n \longrightarrow \mathbb{C}, \quad \mathscr{Q}_n(\alpha, \beta) = \mathrm{i}^{n} \,  (-1)^{\binom{n}{2}} \, \mathscr{R}(\alpha, \bar\beta), \quad \text{where $\mathrm{i} =\sqrt{-1}$.}
\]
Then, by the hard Lefschetz theorem and the Hodge--Riemann relations for compact K\"ahler manifolds \cite[Section 3.3]{Huybrechts}, each $(\mathscr{M}_n, \mathscr{Q}_n)$ is a complex Lefschetz module of degree $d-n$ over $(A, \mathscr{K}_A)$. 
\end{example}

We now discuss the behavior of the decomposition package with respect to Hodge structures. As before, let $A$ and $\mathscr{K}_A$ be as in Section~\ref{sec:introduction}, and set $A_{\mathbb{C}}=A\otimes_{\mathbb{R}}\mathbb{C}$.  
Let $H=\bigoplus_{r \ge 0} H^r$ be a finite-dimensional graded real vector space with an $A$-action satisfying $A^k H^r \subseteq H^{r+2k}$, and let $\mathscr{R} \colon H \times H \to \mathbb{R}$ be an $A$-invariant graded-symmetric bilinear form:
\[
\mathscr{R}(ax_1,x_2)=\mathscr{R}(x_1,ax_2) \ \ \text{and} \ \ 
\mathscr{R}(x_1,x_2)=(-1)^{st}\mathscr{R}(x_2,x_1) \ \ \text{for all $a \in A$, $x_1 \in H^s$, and $x_2 \in H^t$.}
\]
We write $H_\mathbb{C}$ for the complexification of $H$, and we extend $\mathscr{R}$ to an $A_\mathbb{C}$-invariant graded-symmetric complex bilinear form $H_\mathbb{C} \times H_\mathbb{C} \to \mathbb{C}$.

\begin{definition}
We say that $H$ is a \emph{graded $A$-module with a Hodge structure} if the following holds:
\begin{enumerate}\itemsep 5pt
  \item[(HS1)]\label{item:HS1} For each $r$, the graded piece $H^r_{\mathbb{C}}$ is equipped with an effective pure Hodge structure of weight $r$, and the elements of $A_\mathbb{C}^k$ act as morphisms of Hodge structures of bidegree $(k,k)$:
\[
H_{\mathbb{C}}^r = \bigoplus_{\substack{p+q =r, \\ p, q \ge 0}} H^{p,q}, \ \  \overline{H^{p,q}} = H^{q,p}, \ \  \text{and} \ \  A_\mathbb{C}^{k}\ H^{p,q} \subseteq H^{k+p, k+q}.
\]
\end{enumerate}
For each nonnegative $n$, consider the graded $A_\mathbb{C}$-modules 
\[
\mathscr{H}_n = \bigoplus_{k\geq 0} \mathscr{H}_n^k \ \ \text{with} \ \ \mathscr{H}_n^k \coloneq H^{k+n,k} \ \ \text{and} \ \ 
\mathscr{H}_{-n} = \bigoplus_{k\geq 0} \mathscr{H}_{-n}^k \ \ \text{with} \ \ \mathscr{H}_{-n}^k \coloneq H^{k,k+n}. 
\]
We say that $(H, \mathscr{R})$ is a \emph{Lefschetz module with a Hodge structure} of degree $d$ over $(A,\mathscr{K}_A)$ if 
the following holds:
\begin{enumerate}
\item[(HS2)] 
The induced bilinear pairing between the $A_\mathbb{C}$-submodules
\[
\mathscr{R} \colon \mathscr{H}_i \times \mathscr{H}_j \longrightarrow \mathbb{C}
\]
is nondegenerate if $i+j=0$, and it is zero if $i+j \neq 0$.
\item[(HS3)]\label{item:HS2} For each nonnegative $n$, consider the Hermitian form 
\[
\mathscr{Q}_n \colon \mathscr{H}_n \times \mathscr{H}_n \longrightarrow \mathbb{C}, \quad \mathscr{Q}_n(\alpha, \beta) = \mathrm{i}^{n} \,  (-1)^{\binom{n}{2}} \,  \mathscr{R}( \alpha, \bar\beta ).
\]
Then $(\mathscr{H}_n, \mathscr{Q}_n)$ is a complex Lefschetz module of degree $d-|n|$ over $(A, \mathscr{K}_A)$.  
\end{enumerate}
\end{definition}

The category of finite-dimensional graded $A$-modules with Hodge structures is an abelian category in which every object has finite length, so the Krull--Schmidt theorem applies to it. 
Since $\mathscr{H}_0$  and $\mathscr{H}_n \oplus \mathscr{H}_{-n}$ for positive $n$ are defined over $\mathbb{R}$, there are graded $A$-submodules $H_0$ and $H_n$ of $H$ such that $H_0 \otimes_\mathbb{R} \mathbb{C} =\mathscr{H}_0$ and 
 $H_n \otimes_\mathbb{R} \mathbb{C} =\mathscr{H}_n \oplus \mathscr{H}_{-n}$.\footnote{In our current convention, the degree $2k+n$ part of $H_n \otimes_\mathbb{R} \mathbb{C}$ corresponds to the degree $k$ part of $\mathscr{H}_n \oplus \mathscr{H}_{-n}$.} 
 Each $H_n$ is a graded $A$-module with a Hodge structure, and we have a decomposition in the category of graded $A$-modules with Hodge structures
 \[
 H = \bigoplus_{n\geq 0} H_n.
 \]
 The object $H_0$ is indecomposable in the  category of graded $A$-modules with Hodge structures if and only if $H_0$ is an indecomposable graded $A$-module. For positive $n$,  the object $H_n$ is indecomposable in the  category of graded $A$-modules with Hodge structures if and only if $\mathscr{H}_n$ is an indecomposable graded $A_\mathbb{C}$-module.
Note that $(H,\mathscr{R})$ is a Lefschetz module with a Hodge structure of degree $d$ over $(A,\mathscr{K}_A)$ if and only if the same is true for $(H_n,\mathscr{R})$ for all $n \ge 0$.
 If $H=H_0$, we say that $H$ is of \emph{Hodge--Tate type}.

Let $B$ and $\mathscr K_B$ be as in Section~\ref{sec:introduction}. 
Choose any decomposition of $H$ in the category of graded $B$-modules with Hodge structures, and let $N_\alpha$ be an indecomposable summand of $H$. 
We will construct a $B$-invariant bilinear form $\mathscr{R}_{\alpha}$ on $N_{\alpha}$ for each $\alpha$.
Either $N_{\alpha}$ is of Hodge--Tate type, or $N_{\alpha}$ is a summand of $H_n$ for some $n > 0$. If $N_{\alpha}$ is of Hodge--Tate type, then Theorem~\ref{thmPD} shows that, up to a nonzero real multiple, there is a unique $B$-invariant symmetric bilinear form $\mathscr{Q}_{\alpha}$ which satisfies the degree orthogonality relation; set $\mathscr{R}_{\alpha}$ to be this form. If $N_{\alpha}$ is a summand of $H_n$ for some $n > 0$, then the complexification $N_{\alpha, \mathbb{C}}$ splits as $(N_{\alpha, \mathbb{C}} \cap \mathscr{H}_n) \oplus (N_{\alpha, \mathbb{C}} \cap \mathscr{H}_{-n})$. The graded $B_{\mathbb{C}}$-module $N_{\alpha, \mathbb{C}} \cap \mathscr{H}_n$ is an indecomposable summand of the complex Lefschetz module $\mathscr{H}_n$, and so Theorem~\ref{thmHomcomplex} equips it with a Hermitian form $\mathscr{Q}_{\alpha}$ which gives it the structure of a complex Lefschetz module. Define the form $\mathscr{R}_{\alpha}$ on $N_{\alpha, \mathbb{C}}$ by the following formulas: for $v, w \in N_{\alpha, \mathbb{C}} \cap \mathscr{H}_n$, set
\[
\mathscr{R}_{\alpha}(v, w) = \mathscr{R}_{\alpha}(\overline{v}, \overline{w}) = 0, \, \mathscr{R}_{\alpha}(v, \overline{w}) = i^{-n} (-1)^{\binom{n}{2}} \mathscr{Q}_{\alpha}(v,w), \, \text{and }\mathscr{R}_{\alpha}(\overline{v}, w) = \overline{\mathscr{R}_{\alpha}(v, \overline{w})}.
\]
This form is the complexification of a real-valued $B$-invariant form on $N_{\alpha}$, which is symmetric if $n$ is even and is skew-symmetric if $n$ is odd. 

Choose some $\ell \in \mathscr{K}_B$. Then we can define the perverse filtration of $H$ as before using a graded $\mathbb{R}[\ell]$-module decomposition. 
Let $\operatorname{Gr}$ be the associated graded of $H$ with respect to the perverse filtration.

\begin{theorem}\label{thmDecompositionHS} 
 For each $\alpha$, there is a unique sign $\epsilon_{\alpha}=\pm 1$ such that $(N_{\alpha}, \epsilon_{\alpha} \mathscr{R}_{\alpha})$ is a Lefschetz module with Hodge structure over $(B, \mathscr{K}_B)$.
Moreover, the perverse filtration of $H$ over $B$ is by Hodge substructures, and the graded pieces of the perverse filtration inherit pure Hodge structures. For each $\eta \in \mathscr{K}_{A/B}$ and $\ell \in \mathscr{K}_{B}$, the primitive pieces of $\operatorname{Gr}$ carry polarized Hodge structures. 
\end{theorem}
\begin{proof} 
Each $N_{\alpha}$ is a direct summand of some $H_n$. As described above when we were constructing the bilinear form $\mathscr{R}_{\alpha}$, the first statement follows from Theorems~\ref{thmPD}, \ref{thmHL}, and \ref{thmHR} if $n = 0$, and it follows from their complex analogues if $n > 0$. To see that the perverse filtration is by Hodge substructures, note that, because each $N_{\alpha}$ satisfies the hard Lefschetz property, we have
$$P_j = \bigoplus_{d(\alpha) + k \le j} N_{\alpha}[-k]^{\oplus m(\alpha, k)}.$$
As the degree $i$ piece of $N_{\alpha}$ has a pure Hodge structure, this shows that $\operatorname{Gr}^{i,j}$ inherits a pure Hodge structure. 

The last part of the theorem amounts to checking the appropriate versions of the relative hard Lefschetz property and the relative Hodge--Riemann relations. It suffices to check these properties for each $H_n$. When $n=0$, the relative hard Lefschetz property and the relative Hodge--Riemann relations follow from Theorems~\ref{thmRHL} and \ref{thmRHR}; when $n > 0$, they follow from their complex analogues. 
\end{proof}

\subsubsection{Hodge theory of the decomposition theorem} 
Let $Y$ be a complex projective variety. For a local system $\mathcal{L}$ on a smooth connected locally closed subset $Z$ of $Y$, let $IC(\overline{Z}, \mathcal{L})$ denote the intersection cohomology perverse sheaf in  the bounded derived category of constructible sheaves $D^b_c(Y)$. If $\mathcal{L}$ is a simple local system, then $IC(\overline{Z}, \mathcal{L})$ is a simple object in the category of perverse sheaves on $Y$.  
Let $\operatorname{IH}(\overline{Z}, \mathcal{L}) = \mathcal{H}(IC(\overline{Z}, \mathcal{L})[-\dim Z])$ be the hypercohomology of $IC(\overline{Z}, \mathcal{L})$, shifted so that its nonzero degrees lie between $0$ and $2 \dim Z$.

Let $H(Y)$ denote the subalgebra of $\operatorname{H}^{\bullet}(Y; \mathbb{R})$ generated by the Chern classes of line bundles on $Y$. Then each $\operatorname{IH}(\overline{Z}, \mathcal{L})$ is a graded module over $H(Y)$. Let $\mathscr{K}_Y$ denote the open convex cone in $H^2(Y)$ generated by the first Chern classes of ample line bundles. 

Now let $f \colon X \to Y$ be a map from $X$ to a projective variety $Y$, and assume that $X$ is smooth of dimension $d$. The decomposition theorem of Beilinson, Bernstein, Deligne, and Gabber \cite[Theorem 6.2.5]{BBD} states that there is an isomorphism in $D^b_c(Y)$
$$f_* \underline{\mathbb{R}}_X \simeq \bigoplus_i IC(\overline{Z_i}, \mathcal{L}_i)[-d-e_i],$$
where $\underline{\mathbb{R}}_X$ is the constant sheaf on $X$, the $Z_i$ are smooth connected locally closed subvarieties of $Y$ of dimension $d_i$, each $\mathcal{L}_i$ is a simple local system on $Z_i$, and the $e_i$ are integers. Taking hypercohomology, for each $k$ we obtain a direct sum decomposition
\begin{equation}\label{eq:decompIH}
\operatorname{H}^k(X; \mathbb{R}) \simeq \bigoplus_i \operatorname{IH}^{k - d - e_i + d_i}(\overline{Z_i}, \mathcal{L}_i).
\end{equation}
Both $\operatorname{H}(X; \mathbb{R})$ and $\operatorname{IH}(\overline{Z_i}, \mathcal{L}_i)$ are graded $H(Y)$-modules, and the above isomorphism is an isomorphism of graded $H(Y)$-modules. 

Because $X$ is smooth, the main theorems of Hodge theory give $\operatorname{H}^{\bullet}(X; \mathbb{R})$ the structure of a Lefschetz module with Hodge structure over $H(X)$; see, for example, \cite[Chapter 3]{Huybrechts}. We will consider the decomposition of $\operatorname{H}^{\bullet}(X; \mathbb{R})$ as a module over the image of $H(Y)$ in $H(X)$. 

The perverse truncation functors on $D^b_c(Y)$ induce a filtration $P_{\bullet}'$ on $\operatorname{H}^k(X; \mathbb{R})$, which we call the geometric perverse filtration. In terms of the chosen decomposition \eqref{eq:decompIH}, this filtration is given by
\[
P_j' = \bigoplus_{e_i + d \leq j} \operatorname{IH}(\overline{Z_i}, \mathcal{L}_i).
\]
We also have a perverse filtration $P_{\bullet}$ on the Lefschetz module $\operatorname{H}^{\bullet}(X; \mathbb{R})$, obtained by decomposing $\operatorname{H}^{\bullet}(X; \mathbb{R})$ in the category of graded $H(Y)$-modules with Hodge structures. 

\begin{proposition}\label{prop:geometricperverse}
The geometric perverse filtration is equal to the perverse filtration on $\operatorname{H}^{\bullet}(X; \mathbb{R})$. 
\end{proposition}

\begin{proof}
Let $\ell \in H^2(Y)$ be an ample class. Since $\ell$ restricts to an ample class on the closure of each $Z_i$,  \cite[Theorem 6.2.10]{BBD} implies that multiplication by $\ell^{d_i - k}$ induces an isomorphism from $\operatorname{IH}^{k}(\overline{Z_i}, \mathcal{L}_i)$ to $\operatorname{IH}^{2 d_i - k }(\overline{Z_i}, \mathcal{L}_i)$ for each nonnegative $k \le d_i$ and each $i$. In particular, refining the decomposition \eqref{eq:decompIH} into a decomposition of graded $\mathbb{R}[\ell]$-modules, we see that the summand corresponding to $IC(\overline{Z_i}, \mathcal{L}_i)[-d - e_i]$ is contained in $P_{d + e_i}$, and it intersects $P_{d + e_i-1}$ trivially. 
\end{proof}

In particular, this implies that the geometric perverse filtration on $\operatorname{H}^{\bullet}(X; \mathbb{R})$ is ``Hodge-theoretic.'' More precisely, Theorem~\ref{thmDecompositionHS} gives the following. 

\begin{corollary}\label{cor:polarized}
The associated graded pieces of the geometric perverse filtration are endowed with pure Hodge structures, and for any ample classes $\eta \in H(X)$ and $\ell \in H(Y)$, the primitive pieces $\operatorname{Prim}^{i,j}_{\eta, \ell}  \coloneq
\operatorname{ker}(\eta^{d-j+1} \colon \operatorname{Gr}^{i,j} \rightarrow \operatorname{Gr}^{i+2d-2j+2,2d-j+2}) \cap \, \operatorname{ker}(\ell^{j-i+1}\colon \operatorname{Gr}^{i,j} \rightarrow \operatorname{Gr}^{2j-i+2,j})$ carry polarized Hodge structures. 
\end{corollary}

Corollary~\ref{cor:polarized} was first proved by Saito \cite{Saito88, Saito90} using his theory of mixed Hodge modules, and a simpler proof was given by de Cataldo and Migliorini \cite{dCM05}. The argument above deduces Corollary~\ref{cor:polarized} as a formal consequence of the results of \cite{BBD}.

We now explain how the main theorems can be used to put a polarized pure Hodge structure on the intersection cohomology of a complex projective variety. This was first accomplished by Saito \cite{Saito88, Saito90}. Another proof was given by de Cataldo and Migliorini \cite{dCM05}.

Let $Y$ be a complex projective variety of dimension $d$, and let
$\operatorname{IH}(Y)$ denote the intersection cohomology of $Y$ with
coefficients in the trivial local system. Let
\[
f\colon X\longrightarrow Y
\]
be a projective resolution of singularities. By the decomposition theorem
\cite[Theorem 6.2.5]{BBD}, the intersection complex $IC(Y)$ is a direct
summand of
$f_*\underline{\mathbb{R}}_X[d]$. 
Consider the decomposition of the Lefschetz module with Hodge structure
$\operatorname{H}^{\bullet}(X;\mathbb{R})$ over $\operatorname{H}(Y)$.
By Proposition~\ref{prop:geometricperverse},
\[
\operatorname{IH}(Y)\subseteq P_d
\ \ \text{and}\ \ 
\operatorname{IH}(Y)\cap P_{d-1}=0.
\]
Thus, $\operatorname{IH}(Y)$ may be identified with a subspace of
$\operatorname{Gr}^{\bullet,d}$. This identification is canonical because
the perverse cohomology sheaf 
${}^{\mathfrak{p}}\mathcal{H}^0
\bigl(f_*\underline{\mathbb{R}}_X[d]\bigr)$ 
is semisimple.

Let $\eta$ be an ample class on $X$. Multiplication by $\eta$ induces a
morphism
\[
{}^{\mathfrak{p}}\mathcal{H}^{0}
\bigl(f_*\underline{\mathbb{R}}_X[d]\bigr)
\longrightarrow
{}^{\mathfrak{p}}\mathcal{H}^{2}
\bigl(f_*\underline{\mathbb{R}}_X[d]\bigr),
\]
whose kernel contains $IC(Y)$. 
By Corollary~\ref{cor:polarized}, this realizes $\operatorname{IH}(Y)$ as a
summand of a Lefschetz module with Hodge structure over
$\operatorname{H}(Y)$. Moreover, a standard argument shows that
$\operatorname{IH}(Y)$ is a Hodge substructure of
$\operatorname{Gr}^{\bullet,d}$; see
\cite[Proof of Theorem 2.2.1]{dCM05}. The next result therefore follows from the
observation that, if a graded $B$-module summand of a Lefschetz module with
Hodge structure is also a Hodge substructure, then it is 
a Lefschetz module with Hodge structure.

\begin{corollary}
There is a pure Hodge structure on $\operatorname{IH}(Y)$ so that it is a Lefschetz module with Hodge structure over $(H(Y), \mathscr{K}_Y)$.
\end{corollary}

A standard argument shows  that this Hodge structure on $\operatorname{IH}(Y)$ is independent of the choice of resolution $X$. The proof of \cite[Theorem 4.3.5]{dCLefschetzII} shows that it agrees with the Hodge structure  on $\operatorname{IH}(Y)$ constructed by Saito. 

\begin{remark}
The results on the Hodge-theoretic nature of the decomposition theorem \cite{Saito88,Saito90,dCM05} can be used to prove most of the main theorems in the case when $H$ is the intersection cohomology of a $d$-dimensional complex projective variety $X$ which is equipped with a morphism  to a complex projective variety $Y$, $A$ is the subalgebra of $\operatorname{H}^{\bullet}(X; \mathbb{R})$ generated by Chern classes of ample line bundles, $\mathscr{R}$ is the natural pairing on $H$, $\mathscr{K}_A$ is the ample cone of $X$, and $B$ is the subalgebra of $A$ generated by pullbacks of nef line bundles from $Y$. Applying the decomposition theorem as in \eqref{eq:decompIH}, we obtain a decomposition of $H$ into a direct sum of graded $B$-modules that satisfies Theorem~\ref{thmDecompositionHS}; see \cite[Section 2]{dCM05}. These properties are preserved by refining this decomposition into a direct sum of indecomposable graded $B$-modules with Hodge structures. 

On the other hand, if $B$ is the subalgebra of $A$ generated by a face of the nef cone which does not contain the class of a globally generated line bundle in its relative interior, then the conclusions of our main theorems are new and cannot be deduced from the known geometric results in an obvious way. 
\end{remark}

\section{The perverse filtration}\label{sec:perverse}

Let $(M,\mathscr{Q})$ be a Lefschetz module of degree $d$ over $(A,\mathscr{K}_A)$, and let $B$ be an $\mathbb{R}$-subalgebra of $A$   generated by a subset of the closure of $\mathscr{K}_A$ in $A^1$. 
In this section, we establish some fundamental properties of the perverse filtration $P$ of $M$ over $B$.\footnote{As noted in the introduction, Theorem~\ref{thmHL} implies that the perverse filtration can be described using any decomposition of $M$ into indecomposable graded $B$-modules, as in  \eqref{eq:perverse-filtration}. However, we will not be able to prove this until after the main theorems have been established.}

Choose any $\ell \in \mathscr{K}_B$, 
and choose any decomposition of $M$ into indecomposable $\mathbb{R}[\ell]$-modules
\begin{equation}\label{eq:decomposition}
M \simeq \bigoplus_{e} \bigoplus_k C_{e}[-k]^{\oplus m(e, k)},
\end{equation}
where $C_{e}$ is the $(e+1)$-dimensional cyclic graded $\mathbb{R}[\ell]$-module generated by $C_e^0$. 
By definition, $P_j$ is the subspace of $M$ spanned by the summands $C_e[-k]^{\oplus m(e,k)}$ satisfying $e + 2k \le j$. 
From the alternative description of the perverse filtration in   \eqref{eq:perverse-filtration-annihilators}, we see that, for fixed $\ell$, the perverse filtration is independent of the choice of the decomposition \eqref{eq:decomposition}.

\subsection{Module structures}
The following lemma endows the associated graded $\operatorname{Gr}$ of the perverse filtration with an $A$-module structure $*\colon A\times \operatorname{Gr} \to \operatorname{Gr}$ satisfying $A^s*\operatorname{Gr}^{i,j} \subseteq \operatorname{Gr}^{i + s, j + 2s}$. 

\begin{lemma}\label{lem:Amod}
For each $s$ and $j$, we have $A^s P_j \subseteq P_{j + 2s}$. 
\end{lemma}

\begin{proof}
Recall from \eqref{eq:perverse-filtration-annihilators} that the degree $k$ part of $P_j$ is spanned by elements of the form $\ell^{k -c} m$, where $m$ is a degree $c$ element satisfying 
\[
\ell^{j -c + 1-k}\cdot \ell^{k-c} m = \ell^{j-2c+1} \cdot m=0.
\]
Observe that, for any $a \in A^s$, the element 
\[
a \cdot \ell^{k -c} m = \ell^{(s + k) -(s+c)} \cdot a m \ \ \text{is annihilated by $\ell^{j-c+1-k}=\ell^{j + 2s + 1 - (s + k) - (s+c)}$.}
\]
Thus, the product $a \cdot \ell^{k -c} m$ lies in $P_{j + 2s}$ by \eqref{eq:perverse-filtration-annihilators}, as required.
\end{proof}

It is clear that $\ell P_j \subseteq P_j$ for each $j$. 
We show that, in fact,  $P_j$ 
is a $B$-module for each $j$.
The main ingredient is a theorem of Cattani and Kaplan on nilpotent orbits \cite[Theorem 3.3]{CattaniKaplan}.

\begin{proposition} \label{propIndependencePerverseFiltraion}
The perverse filtration $P$ is independent of the choice of $\ell \in \mathscr{K}_B$. 
\end{proposition}

\begin{proof} 
Choose an increasing family of open polyhedral cones $\mathscr{K}_{B}^1 \subseteq \mathscr{K}_{B}^2 \subseteq \cdots$ whose union is $\mathscr{K}_B$.
Choose any $\eta \in \mathscr{K}_A$, and let $\mathscr{K}_{A}^n$ be the open polyhedral cone  generated by $\eta$ and $\mathscr{K}_{B}^n$ for each $n$. 

Recall that the \emph{weight filtration} $W=\big(W_{-d} \subseteq  \cdots \subseteq W_d\big)$ of $M$ associated to $\ell \in \mathscr{K}_B$ is the unique filtration on $M$ satisfying  
$\ell W_i \subseteq W_{i-2}$ for all $i$ and, for all nonnegative $i$, the map
\[
W_i/W_{i-1} \longrightarrow W_{-i}/W_{-i-1}, \quad x \longmapsto \ell^i x
\]
is an isomorphism of vector spaces. 
By \cite[Theorem 3.1]{Cat08}, for each $n$, we may apply \cite[Theorem 3.3]{CattaniKaplan} to the Lefschetz module $(M,\mathscr{Q})$ over $(A,\mathscr{K}_{A}^n)$. Since  $\mathscr{K}_B=\bigcup_n \mathscr{K}_{B}^n$, the conclusion is  that  $W$ is independent  of  $\ell \in \mathscr{K}_B$. 
The weight filtration determines the perverse filtration by the formula
$P_j = \sum_{i+2k \leq j} W_i \cap M^k$, and hence $P$ is 
 independent of $\ell \in \mathscr{K}_B$. 
\end{proof}

Since $\mathscr{K}_B$ generates $B$ as an $\mathbb{R}$-algebra, 
Proposition~\ref{propIndependencePerverseFiltraion} implies that the action of $B$ preserves the perverse filtration $P$. 

\begin{corollary}\label{cor:B-module}
For any $b \in B$ and any $j$, we have $b P_j \subseteq P_j$. 
\end{corollary}

Corollary~\ref{cor:B-module} implies that each $\operatorname{Gr}^{\bullet, j}$ has the structure of a graded $B$-module, and it satisfies  the following hard Lefschetz property for $\mathscr{K}_B$.

\begin{corollary}\label{cor:lefschetzforell}
For every nonnegative $2i \le j \le 2d$ and $\ell \in \mathscr{K}_B$, the linear map
\[
\operatorname{Gr}^{i,j} \longrightarrow \operatorname{Gr}^{j - i, j}, \qquad x \longmapsto \ell^{j-2i} x
\]
is an isomorphism of vector spaces. 
\end{corollary}

\subsection{The descent lemma}\label{ssec:descent}

For each $a \in A^k$,  consider the $A$-module $M_a \coloneq M/\operatorname{ann}(a)$ and the quotient map $\varphi_a \colon M \to M_a$. The $A$-module $M_a$ is equipped with an $A$-invariant symmetric bilinear form ${\mathscr{Q}}_a \colon {M}_a \times {M}_a \to \mathbb{R}$ defined by 
\[
{\mathscr{Q}}_a(\varphi_a(x_1), \varphi_a(x_2)) = \mathscr{Q}(x_1, ax_2)\ \  \text{for all $x_1, x_2\in M$}.
\]
By construction, the induced bilinear pairing between the graded pieces
\[
{\mathscr{Q}}_a \colon {M}_a^i \times {M}_a^j \longrightarrow \mathbb{R}
\]
is nondegenerate if $i+j=d-k$ and zero if $i+j \neq d-k$.
The next result was proved in \cite[Lemma 1.16]{CKS87} and \cite[Theorem 3.2]{Cat08}, where it is called the \emph{descent lemma}, and also in \cite[Theorem 2.1.5]{KK87}, where it is called the \emph{vanishing cycle theorem}. 

\begin{lemma}\label{lem:descent}
If $\gamma$ is an element of the closure of $\mathscr{K}_A$, then  $(M_{\gamma}, {\mathscr{Q}}_{\gamma})$  is a Lefschetz module of degree $d-1$ over $(A, \mathscr{K}_A)$. 
\end{lemma}

We will frequently consider $M_{\ell^p} = M/\operatorname{ann}(\ell^p)$ for a nonnegative integer $p$. By repeatedly applying Lemma~\ref{lem:descent}, we deduce that $(M_{\ell^p},\mathscr{Q}_{\ell^p})$ is a Lefschetz module of degree $d-p$ over $(A, \mathscr{K}_A)$. 
The chosen decomposition~\eqref{eq:decomposition} of $M$ induces the following decomposition of $M_{\ell^p}$ into indecomposable $\mathbb{R}[\ell]$-modules: 
\begin{equation}\label{eq:decomposition-descent}
M_{\ell^p} \simeq \bigoplus_{e \geq p} \bigoplus_k C_{e-p}[-k]^{\oplus m(e, k)}.
\end{equation}

\subsection{The 1-dimensional summands} \label{sec:submodule-N}
Let $N$ 
be the direct sum of the $1$-dimensional $\mathbb{R}[\ell]$-module summands in the chosen decomposition \eqref{eq:decomposition} of $M$. The restriction of the perverse filtration of $M$ over $B$ to $N$ coincides with the filtration induced by the grading of $M$, and hence $N$ can be identified with its image in $\operatorname{Gr}$. 
The image of $N$ does not depend on the decomposition \eqref{eq:decomposition}, and it admits the description
\[
N \simeq \operatorname{ann}_M(\ell) \big/ \bigl(\ell M \cap \operatorname{ann}_M(\ell)\bigr),
\]
where $\operatorname{ann}_M(\ell)$ is the kernel of multiplication by $\ell$ in $M$. 
In fact, the image of $N$ coincides with the kernel of the map $\operatorname{Gr}^{\bullet, 2\bullet} \to \operatorname{Gr}^{\bullet + 1, 2 \bullet}$ induced by multiplication by $\ell$, and it is therefore an $A$-submodule of $\operatorname{Gr}$.

\subsection{The bilinear form $\underline{\mathscr{Q}}$ on $\operatorname{Gr}$}

We check that $\mathscr{Q}$ descends to a symmetric bilinear form 
\[
\underline{\mathscr{Q}} \colon \operatorname{Gr}^{i,j} \times \operatorname{Gr}^{d - i, 2d - j} \longrightarrow \mathbb{R}.
\]

\begin{lemma}\label{lem:formdefined}
For all $i$ and $j$, we have
$\mathscr{Q}(M^i \cap P_j, M^{d - i} \cap P_{2d - j - 1}) =0$.
\end{lemma}
\begin{proof}
Note that $M^i \cap P_j$ is spanned by elements of the form $\ell^{i-c}  m$, where $m \in M^c$ with $\ell^{e+1}  m = 0$ for some $e$ with $e + 2c \le j$. Similarly, $M^{d-i} \cap P_{2d - j - 1}$ is spanned by elements of the form $\ell^{d - i - c'}  m'$ with $\ell^{e'+1}  m' = 0$ for some $e'$ with $e' + 2c' \le 2d - j - 1$. Then
\[
\mathscr{Q}(\ell^{i - c}  m, \ell^{d - i - c'} m') = \mathscr{Q}(\ell^{d - c - c'}m, m') = \mathscr{Q}(m, \ell^{d - c - c'}m').
\]
This vanishes if either $e +1 \le d- c - c' $ or $e' + 1 \le d-c -c'$. Adding the two inequalities $e + 2c \le j$ and $e' + 2c' \le 2d - j - 1$ implies that one of these must hold.  
\end{proof}

We extend $\underline{\mathscr{Q}}$ to a map $\operatorname{Gr} \times \operatorname{Gr} \to \mathbb{R}$ by setting $\underline{\mathscr{Q}}(\operatorname{Gr}^{i,j} \times \operatorname{Gr}^{i', j'}) = 0$ unless $i + i' = d$ and $j + j' = 2d$.

Recall from Section~\ref{sec:submodule-N} that $N$ is the subspace of $M$ spanned by the $1$-dimensional summands in the chosen $\mathbb{R}[\ell]$-module decomposition, and that it can be identified with $\operatorname{ann}_M(\ell)/(\ell M \cap \operatorname{ann}_M(\ell))$, which is naturally a subspace of $\operatorname{Gr}$. 
Since the perverse filtration on $N$ coincides with the filtration induced by the grading, the restrictions of the forms $\mathscr{Q}$ and $\underline{\mathscr{Q}}$ to $N$ can be identified. 

\begin{lemma}\label{lem:nondegenN}
The restriction of $\mathscr{Q}$ to $N$ is nondegenerate. 
\end{lemma}

\begin{proof}
We need to show that each element $x \in \operatorname{ann}_M(\ell)$ which has $\mathscr{Q}(x, y) = 0$ for all $y \in \operatorname{ann}_M(\ell)$ lies in $\ell M$. Note that $\operatorname{ann}_M(\ell)$ is the orthogonal complement of $\ell M$ with respect to $\mathscr{Q}$. Then  Poincar\'e duality for $\mathscr{Q}$  
implies that $\ell M$ is the orthogonal complement of $\operatorname{ann}_M(\ell)$.
\end{proof}

In particular, we deduce that $\dim N^i = \dim N^{d-i}$. Recall that $m(e, k)$ is the multiplicity of $C_e[-k]$ in the $\mathbb{R}[\ell]$-module decomposition of $M$. 

\begin{proposition}\label{prop:nondegenerate}
The form $\underline{\mathscr{Q}}$ is nondegenerate. 
\end{proposition}

\begin{proof}
We first show that $\dim \operatorname{Gr}^{i,j} = \dim \operatorname{Gr}^{d - i, 2d-j}$ for all $i$ and $j$. It follows from Lemma~\ref{lem:nondegenN} that $m(0, k) = m(0, d-k)$. By \eqref{eq:decomposition-descent}, $m(e, k)$ is equal to the multiplicity of $C_0[-k]$ in $M_{\ell^e}$. Using Lemma~\ref{lem:descent} and applying Lemma~\ref{lem:nondegenN} to $M_{\ell^e}$, we deduce that $m(e, k) = m(e, d - e - k)$. This implies that $\dim \operatorname{Gr}^{i,j} = \dim \operatorname{Gr}^{d - i, 2d-j}$.

For each $k$, choose a basis for $M^k$ which is compatible with the perverse filtration on $M^k$. Lemma~\ref{lem:formdefined} and the above equality of dimensions imply that the matrix representing $\mathscr{Q}$ with respect to this basis is triangular. As $\mathscr{Q}$ is nondegenerate, this implies that $\underline{\mathscr{Q}}$ is nondegenerate. 
\end{proof}

\subsection{Applications of the descent lemma}
For $\gamma \in \overline{\mathscr{K}}_A$, let $P_{\gamma, \bullet}$ be the perverse filtration on $M_{\gamma}$, defined using our chosen element $\ell \in \mathscr{K}_B$. Let $\operatorname{Gr}^{}_{\gamma}$ be the associated graded of $M_{\gamma}$ with respect to the perverse filtration. Let $\psi \colon M \to M_{\ell}$ be the quotient map associated to $\ell$.

\begin{proposition}\label{prop:ellfunctoriality}
The quotient map $\psi \colon M \to M_{\ell}$ satisfies $\psi(P_j) = P_{\ell, j-1}$ for all $j$. 
 The induced map $\operatorname{Gr} \psi \colon \operatorname{Gr}^{i, j} \to \operatorname{Gr}^{i, j - 1}_{\ell}$ is an isomorphism if $i < j/2$. 
\end{proposition}

\begin{proof}
If $M \simeq \bigoplus_{e, k} C_e[-k]^{\oplus m(e, k)}$ is the fixed $\mathbb{R}[\ell]$-module decomposition of $M$, then $\psi(C_e[-k]) = C_{e-1}[-k]$, where $C_{-1}$ is interpreted as $0$. By \eqref{eq:decomposition-descent}, we have an $\mathbb{R}[\ell]$-module decomposition
\[
M_{\ell} \simeq \bigoplus_{e, k} C_{e-1}[-k]^{\oplus m(e, k)}.
\] 
As $C_e[-k]$ lies in $P_{e + 2k}$ and $C_{e-1}[-k]$ lies in $P_{\ell, e + 2k - 1}$, this implies that $\psi(P_j) = P_{\ell, j-1}$. 

In terms of the chosen decomposition, we have
\[
\operatorname{Gr}^{\bullet, j} = \bigoplus_{e + 2k = j} C_e[-k]^{\oplus m(e,k) } \ \ \text{and} \ \ \operatorname{Gr}^{\bullet, j-1}_{\ell} = \bigoplus_{e + 2k = j} \psi(C_e[-k]^{\oplus m(e, k)}).
\]
The restriction of $\psi$ to $C_e[-k]$ is an isomorphism in degree $i$ if $i < e + k$. Since $j/2  = e/2 +k$ and $i < j/2$, we see that $i < e + k$, giving the result.
\end{proof}

The definitions of the forms $\underline{\mathscr{Q}}$ on $\operatorname{Gr}$ and $\underline{\mathscr{Q}}_{\ell}$ on $\operatorname{Gr}_{\ell}$ in terms of liftings to $M$ and $M_{\ell}$ give the following compatibility.

\begin{corollary}\label{cor:gysinell}
For all $x_1 \in \operatorname{Gr}^{i,j}$ and $x_2 \in \operatorname{Gr}^{d - i - 1, 2d - j}$, we have
\[
\underline{\mathscr{Q}}_{\ell}(\operatorname{Gr} \psi(x_1), \operatorname{Gr} \psi(x_2)) = \underline{\mathscr{Q}}(x_1, \ell x_2).
\]
\end{corollary}

Choose some $\eta \in \mathscr{K}_A$, and let $M_{\eta} = M/\operatorname{ann}(\eta)$. Let $\varphi \colon M \to M_{\eta}$ be the quotient map. Note that the hard Lefschetz property 
implies that $\varphi$ is an isomorphism in every degree less than $d/2$. 

\begin{proposition}\label{prop:etagysin}
The map $\varphi \colon M \to M_{\eta}$ satisfies $\varphi(P_j) \subseteq P_{\eta,j}$, and the induced map $\operatorname{Gr} \varphi \colon \operatorname{Gr}^{i,j} \to \operatorname{Gr}^{i,j}_{\eta}$ is injective for $j < d$. 
\end{proposition}

\begin{proof}
That $\varphi(P_j) \subseteq P_{\eta,j}$ follows from the description of the perverse filtration in \eqref{eq:perverse-filtration-annihilators}. Choose some $\ell \in \mathscr{K}_B$, and note that Corollary~\ref{cor:lefschetzforell} implies that $\operatorname{Gr}^{\bullet,j}$ and $\operatorname{Gr}^{\bullet, j}_{\eta}$ are representations of $\mathfrak{sl}_2$, where $\ell$ acts as a raising operator. Furthermore, the map $\operatorname{Gr} \varphi$ is a map of representations of $\mathfrak{sl}_2$. In particular, in order to check that it is injective for $j < d$, it suffices to check that it is injective when restricted to $\operatorname{ker}(\ell^{j - 2i + 1} \colon \operatorname{Gr}^{i,j} \to \operatorname{Gr}^{j - i + 1, j})$ for each $i \le j/2$. 

If $i < j/2$, then, using Proposition~\ref{prop:ellfunctoriality}, we may replace $M$ by $M/\operatorname{ann}_M(\ell^{j - 2i})$. By Lemma~\ref{lem:descent}, the result is still a Lefschetz module, and hence we may assume that $j = 2i$. As in  Section~\ref{sec:submodule-N}, we have an identification of $A$-modules
\[
\bigoplus_i \ker(\ell \colon  \operatorname{Gr}^{i, 2i} \to \operatorname{Gr}^{i + 1, 2i}) = \operatorname{ann}_M(\ell) \big/ \big( \ell M \cap \operatorname{ann}_M(\ell) \big).
\] 
Suppose we have $x \in \ker(\ell \colon  \operatorname{Gr}^{i, 2i} \to \operatorname{Gr}^{i + 1, 2i})$ with $i < d/2$. Then $\operatorname{Gr} \varphi(x) = 0$ if and only if $x$ is in the kernel of the map
\[
\operatorname{ann}_M(\ell) \big/ \big(\ell M \cap \operatorname{ann}_M(\ell) \big) \longrightarrow \operatorname{ann}_{M_{\eta}}(\ell) \big/ \big(\ell M_{\eta} \cap \operatorname{ann}_{M_{\eta}}(\ell)\big).
\]
In other words, lifting $x$ to $\tilde{x} \in M$, we have $\varphi(\tilde{x}) = \ell \cdot y$ for some $y \in M^{i-1}_{\eta}$. Using Poincar\'e duality, 
we can identify $M_{\eta}[-1]$ with $\eta M \subseteq M$. Then this means that $\eta \cdot \tilde{x} = \eta \cdot \ell \cdot y$ in $M^{i+1}$. 
Since $i < d/2$, by the hard Lefschetz property,  this implies that $\tilde{x} = \ell \cdot y$, so $x$ vanishes in $\operatorname{Gr}^{i, 2i}$, as required.
\end{proof}

From the definitions of the bilinear forms $\underline{\mathscr{Q}}$ and $\underline{\mathscr{Q}}_\eta$, we have the following compatibility. 

\begin{corollary}\label{cor:etagysin}
For all $x_1 \in \operatorname{Gr}^{i,j}$ and $x_2 \in \operatorname{Gr}^{d - i - 1, 2d - j - 2}$, we have
\[
\underline{\mathscr{Q}}_{\eta}(\operatorname{Gr} \varphi(x_1), \operatorname{Gr} \varphi(x_2)) = \underline{\mathscr{Q}}(x_1, \eta * x_2).
\]
\end{corollary}

\section{Relative hard Lefschetz property and relative Hodge--Riemann relations}\label{sec:relative}

In this section, we prove Theorem~\ref{thmRHL} (relative hard Lefschetz property) and  Theorem~\ref{thmRHR} (relative Hodge--Riemann relations).
The two statements are established
 simultaneously by induction on the degree $d$ of the Lefschetz module $(M,\mathscr{Q})$ over $(A,\mathscr{K}_A)$.

\begin{proposition}\label{prop:raised}
Suppose that Theorems~\ref{thmRHL} and~\ref{thmRHR} hold for every Lefschetz module of degree at most $d-1$. Then, for every nonnegative $j\le d$ and  $\eta\in \mathscr{K}_{A}$, the linear map
\[
\operatorname{Gr}^{\bullet,j} \longrightarrow \operatorname{Gr}^{\bullet+d-j,2d-j}, \qquad x \longmapsto \eta^{d-j} *  x 
\]
is an isomorphism of vector spaces.
\end{proposition}

\begin{proof}
We first show that the map is injective. 
Choose some $\ell \in \mathscr{K}_B$. By Corollary~\ref{cor:lefschetzforell}, there are $\mathfrak{sl}_2$ actions on $\operatorname{Gr}^{\bullet,j}$ and $\operatorname{Gr}^{\bullet+d-j,2d-j}$ where $\ell$ acts as the raising operator, and $\eta^{d-j} *$ is a map of representations of $\mathfrak{sl}_2$. In order to check that this map is injective, it suffices to check that it is injective on $\ker(\ell^{j - 2i + 1} \colon \operatorname{Gr}^{i,j} \to \operatorname{Gr}^{j - i + 1, j})$ for each $i\leq j/2$. If $j = d$, then $\eta^{d - j}*$ is the identity, so we may assume that $j < d$. 
Suppose we have some element
\[
x \in \ker( \operatorname{Gr}^{i,j} \xrightarrow{\ell^{j - 2i + 1}} \operatorname{Gr}^{j - i + 1, j}) \cap \ker( \operatorname{Gr}^{i, j} \xrightarrow{\eta^{d-j} *} \operatorname{Gr}^{i + d - j, 2d-j}).
\]
Let $\varphi \colon M \to M_{\eta} = M/\operatorname{ann}_M(\eta)$ be the quotient map, and let $\operatorname{Gr} \varphi \colon \operatorname{Gr} \to \operatorname{Gr}_{\eta}$ be the associated graded map in  Proposition~\ref{prop:etagysin}. By construction, we have
\[
\operatorname{Gr} \varphi(x) \in \ker(\operatorname{Gr}^{i,j}_{\eta} \xrightarrow{\ell^{j - 2i + 1} } \operatorname{Gr}^{j - i + 1, j}_{\eta}) \cap \ker( \operatorname{Gr}_{\eta}^{i, j} \xrightarrow{\eta^{d-j} *} \operatorname{Gr}_{\eta}^{i + d - j, 2d-j}).
\]
That is, $\operatorname{Gr} \varphi(x)$ lies in the $(i, j)$-th graded piece of the primitive part of $\operatorname{Gr}_{\eta}$.  
By induction on the degree $d$, we have 
\[
(-1)^{i} \underline{\mathscr Q}_{\eta}(\operatorname{Gr} \varphi(x), \eta^{d-j-1} * \ell^{j - 2i}\operatorname{Gr}\varphi(x)) \ge 0,
\]
with equality if and only if $\operatorname{Gr} \varphi(x) = 0$. 
By Corollary~\ref{cor:etagysin}, we have
\[
\underline{\mathscr Q}_{\eta}(\operatorname{Gr} \varphi(x), \eta^{d-j-1} * \ell^{j - 2i} \operatorname{Gr}\varphi(x)) = \underline{\mathscr Q}(x, \eta^{d-j}*\ell^{j - 2i}x) = 0,
\]
so we get that $\operatorname{Gr} \varphi(x) = 0$. By the injectivity part of Proposition~\ref{prop:etagysin}, it follows that $x = 0$. 

We have shown that $\eta^{d-j} *$ is injective. Multiplying by $\ell^{j - 2i}$ and using Corollary~\ref{cor:lefschetzforell}, we have an injection
\[
\eta^{d-j} * \ell^{j - 2i} \colon  \operatorname{Gr}^{i,j} \hookrightarrow \operatorname{Gr}^{d -i, 2d - j}.
\]
By Proposition~\ref{prop:nondegenerate}, $\dim \operatorname{Gr}^{i,j}  = \dim \operatorname{Gr}^{d -i, 2d - j}$, so this map is an isomorphism. This proves that $\eta^{d-j} *$ is an isomorphism.  
\end{proof}

Once we know that the conclusion of Proposition~\ref{prop:raised} holds for a Lefschetz module $M$,  we may decompose $\operatorname{Gr}$ into primitive pieces as in \eqref{eq:primitive-decomposition}.

\begin{proposition}\label{prop:easyprimitive}
Suppose that Theorems~\ref{thmRHL} and ~\ref{thmRHR} hold for every Lefschetz module of degree at most $d-1$. If $(i, j) \not= (d/2, d)$, then, for any $\eta \in \mathscr{K}_A$ and any $\ell \in \mathscr{K}_B$, the  symmetric bilinear form on $\operatorname{Gr}^{i,j}$ defined by
\[
(x_1,x_2) \longmapsto (-1)^i \underline{\mathscr{Q}}(x_1, \eta^{d-j} * \ell^{j-2i} x_2)
\]
is positive definite when restricted to $\operatorname{Prim}^{i,j}_{\eta,\ell}$.
\end{proposition}

\begin{proof}
Let $x \in \operatorname{Prim}^{i,j}_{\eta,\ell}$ be a nonzero element. First, suppose that $j < d$. Then, by Proposition~\ref{prop:etagysin}, there is an injective map $ \operatorname{Gr} \varphi \colon \operatorname{Gr}^{i,j} \to \operatorname{Gr}_{\eta}^{i,j}$. 
We claim that $\operatorname{Gr} \varphi(x)$ remains primitive. Because the map $\operatorname{Gr} \varphi$ respects both the $A$-module structure and the $B$-module structure, it is clear that 
\[
\ell^{j - 2i + 1} \operatorname{Gr} \varphi(x) = 0.
\]

We check that $\eta^{d - j}* \operatorname{Gr} \varphi(x) = 0$,
which holds if and only if $\eta^{d-j} * \ell^{j - 2i} \operatorname{Gr} \varphi(x) = 0$ by Corollary~\ref{cor:lefschetzforell}.
Using Proposition~\ref{prop:ellfunctoriality}, we may replace $M$ with $M_{\ell^{j - 2i}}$, and we can therefore assume that $j = 2i$. 
As noted in Section~\ref{sec:submodule-N}, we have an identification
\[
N \coloneq \bigoplus_q \ker(\ell \colon \operatorname{Gr}^{q,2q} \to \operatorname{Gr}^{q+1, 2q}) = \operatorname{ann}_M(\ell) \big/ \big( \ell M \cap \operatorname{ann}_M(\ell) \big).
\]
It then suffices to note that, if $x$ lies in $N^i$ with $\eta^{d -2i + 1} * x = 0$, then the image of $x$ in $N/\operatorname{ann}_N(\eta)$ is annihilated by $\eta^{d - 2i}$. 
Since Theorem~\ref{thmRHR} holds for $M_{\eta}$,  we have 
\[
(-1)^i \underline{\mathscr{Q}}_{\eta}(\operatorname{Gr} \varphi(x), \eta^{d-j-1} * \ell^{j-2i} \operatorname{Gr} \varphi(x)) > 0.
\]
The result then follows from Corollary~\ref{cor:etagysin}. 

Now suppose that $j = d$, but $i < d/2$. By Proposition~\ref{prop:ellfunctoriality}, we have an isomorphism \[
\operatorname{Gr} \psi \colon \operatorname{Gr}^{i,d} \xrightarrow{\simeq} \operatorname{Gr}^{i,d-1}_{\ell}.
\]
Since this map respects both the $A$-module structure and the $B$-module structure, it is clear that $\eta * \operatorname{Gr} \psi(x) = 0$. It follows from the description of the perverse filtration in \eqref{eq:decomposition} that $\ell^{j - 2i} \operatorname{Gr} \psi(x) = 0$, so $\operatorname{Gr} \psi(x)$ is primitive. Since Theorem~\ref{thmRHR} holds for $M_{\ell}$, the result follows from Corollary~\ref{cor:gysinell}.
\end{proof}

We now consider the remaining primitive part  $\operatorname{Prim}^{d/2, d}_{\eta,\ell}$ in the case when $d$ is even. Our approach is to analyze the signature of the restriction of $\underline{\mathscr{Q}}$ to $\operatorname{Gr}^{d/2, \bullet}$. An alternative approach to a similar problem is given in \cite[Section 5.4]{dCM05}. 
We begin by computing the signature of the restriction of $\underline{\mathscr{Q}}$ to $\operatorname{Gr}^{d/2, \bullet}$. 

\begin{proposition}\label{prop:qunderlinesignature}
The signature of the restriction of $\underline{\mathscr{Q}}$ to $\operatorname{Gr}^{d/2, \bullet}$ is equal to  
\[
\sum_{i=0}^{d/2} (-1)^i (\dim M^i - \dim M^{i - 1}).
\]
\end{proposition}

\begin{proof}
Let $b_j = \dim \operatorname{Gr}^{d/2, j}$. 
By Proposition~\ref{prop:nondegenerate}, we have $b_j = b_{2d - j}$. 
Choose a basis $x_1, x_2, \dotsc, x_r$ for $M^{d/2}$ such that 
\[
x_{r_{j-1} + 1}, \dotsc, x_{r_{j}} \in P_j \cap M^{d/2} \ \ \text{for all $j$},
\]
where $r = \dim M^{d/2}$ and  $r_j = \sum_{e\leq j}b_e$.
Then the images of $x_i$ form a basis for $\operatorname{Gr}^{d/2, \bullet}$.

Let $T$ be the symmetric matrix defined by $T_{s,t} = \mathscr{Q}(x_{s}, x_t)$, so that $T$ represents the restriction of $\mathscr{Q}$ to $M^{d/2}$. The Hodge--Riemann relations 
imply that the signature of $T$ is 
\[
\sum_{i=0}^{d/2} (-1)^i (\dim M^i - \dim M^{i - 1}).
\]
We can decompose $T$ into blocks, representing how different pieces of the perverse filtration pair with each other. 
By Lemma~\ref{lem:formdefined}, $T$ vanishes northwest of the blocks along the main antidiagonal, which represent the pairing of complementary pieces of the perverse filtration. 

The matrix representing the restriction of $\underline{\mathscr{Q}}$ to $\operatorname{Gr}^{d/2, \bullet}$ is obtained by setting the entries in the blocks below the antidiagonal to $0$. In particular, $\underline{\mathscr{Q}}$ has the same signature as $\mathscr{Q}$, giving the result. 
\end{proof}

Suppose that $d$ is even and that Theorems~\ref{thmRHL} and ~\ref{thmRHR} hold for all Lefschetz modules of degree at most $d-1$.  Propositions~\ref{prop:raised} and ~\ref{prop:easyprimitive} then give a decomposition
\[
\operatorname{Gr}^{d/2, \bullet} = \bigoplus_{\substack{j \le d \\ i \le j/2}} \bigoplus_{\substack{s \le d - j \\ t \le j - 2i \\ s + t = d/2 - i}} \eta^{s} * \ell^t \operatorname{Prim}^{i,j}_{\eta,\ell}.
\]
This decomposition is nearly orthogonal with respect to $\underline{\mathscr{Q}}$: we have
\[
\underline{\mathscr{Q}}(\eta^{s} * \ell^t \operatorname{Prim}^{i,j}_{\eta,\ell}, \eta^{s'} * \ell^{t'} \operatorname{Prim}^{i',j'}_{\eta,\ell}) = 0,
\]
unless $i = i'$, $j = j'$, $s + s' = d - j$, and $t + t' = j - 2i$. 
We pair the summand $\eta^{s} *\ell^t \operatorname{Prim}^{i,j}_{\eta,\ell}$ with the complementary summand $\eta^{d - j - s} * \ell^{j - 2i - t} \operatorname{Prim}^{i,j}_{\eta,\ell}$. 
Whenever these two summands are distinct, the restriction of
\(\underline{\mathscr{Q}}\) to their direct sum has a matrix of the form
\[
\begin{pmatrix}
0 & U\\
U^{t} & 0
\end{pmatrix},
\]
whose signature is zero. Thus, only the self-complementary summands can contribute to the signature, and 
 the signature of $\underline{\mathscr{Q}}$ is the same as the signature of the restriction of $\underline{\mathscr{Q}}$ to $\operatorname{Gr}^{d/2, d}$. 
The decomposition of $\operatorname{Gr}^{d/2, d}$ takes the simpler form
\[
\operatorname{Gr}^{d/2, d} = \bigoplus_{\substack{j \le d \text{ even} \\  i \le j/2}} \eta^{(d - j)/2}* \ell^{j/2 - i} \operatorname{Prim}^{i,j}_{\eta,\ell}.
\]
The following identity compares the signatures of the two sides of this decomposition.

\begin{lemma}\label{lem:numerical}
We have
\[
\sum_{i=0}^{d/2} (-1)^i (\dim M^i - \dim M^{i - 1}) = \sum_{\substack{j \le d \text{ even} \\ i \le j/2}} (-1)^i \dim \operatorname{Prim}^{i,j}_{\eta,\ell}.
\]
\end{lemma}

\begin{proof}
We decompose $\operatorname{Gr}$ as in \eqref{eq:primitive-decomposition} and compute the contribution to the left-hand side of the terms associated to $\operatorname{Prim}^{i,j}_{\eta,\ell}$. For each $p$, the dimension of the sum of the pieces associated to $\operatorname{Prim}^{i,j}_{\eta,\ell}$ in $\operatorname{Gr}^{p, \bullet}$ is $\dim \operatorname{Prim}^{i,j}_{\eta,\ell}$ times the number of pairs $(s, t)$ with $s + t = p - i$, $0 \le s \le d - j$, and $0 \le t \le j - 2i$. 

The number of such pairs is $0$ for $p < i$, increases by one when we increase $p$ by one until $p$ reaches $i + \operatorname{min}(d - j, j - 2i)$, which happens for some $p \le d/2$, and then is constant until $p$ reaches $d/2$. Therefore, the contribution of the terms associated to $\operatorname{Prim}^{i,j}_{\eta,\ell}$ is
\[
\sum_{0 \le q \le \operatorname{min}(d - j, j- 2i)} (-1)^{i + q} \dim \operatorname{Prim}^{i,j}_{\eta,\ell}.
\]
Since $d$ is even, $d - j$ and $j - 2i$ both have the same parity as $j$, so this sum is $0$ if $j$ is odd and is $(-1)^{i} \dim \operatorname{Prim}^{i,j}_{\eta,\ell}$ if $j$ is even. 
\end{proof}

We continue to assume that $d$ is even and that Theorems~\ref{thmRHL} and ~\ref{thmRHR} hold for every Lefschetz module of degree at most $d-1$.

\begin{proposition}\label{prop:criticalprimitive}
For any $\eta \in \mathscr{K}_A$ and $\ell \in \mathscr{K}_B$, the restriction of $(-1)^{d/2} \, \underline{\mathscr{Q}}$ to $\operatorname{Prim}^{d/2, d}_{\eta,\ell}$ is positive definite. 
\end{proposition}

\begin{proof}
Let $n$ be the signature of $\underline{\mathscr{Q}}$ on $\operatorname{Prim}^{d/2, d}_{\eta,\ell}$. By Proposition~\ref{prop:qunderlinesignature} and Lemma~\ref{lem:numerical}, the signature of $\underline{\mathscr{Q}}$ on $\operatorname{Gr}^{d/2, \bullet}$ is equal to 
\[
\sum_{i=0}^{d/2} (-1)^i (\dim M^i - \dim M^{i - 1}) = \sum_{\substack{j \le d \text{ even}\\ i \le j/2}} (-1)^i \dim \operatorname{Prim}^{i,j}_{\eta,\ell}.
\]
On the other hand, using Proposition~\ref{prop:easyprimitive} and the discussion above about the decomposition being nearly orthogonal with respect to $\underline{\mathscr{Q}}$, the signature of $\underline{\mathscr{Q}}$ on $\operatorname{Gr}^{d/2, \bullet}$ is
\[
n + \sum_{\substack{j \le d \text{ even}\\ i \le j/2  \\ (i, j) \not= (d/2, d)}} (-1)^i \dim \operatorname{Prim}^{i,j}_{\eta,\ell}.
\]
Comparing these expressions shows that $n = (-1)^{d/2} \dim \operatorname{Prim}^{d/2, d}_{\eta,\ell}$. 
\end{proof}

\begin{proof}[Proofs of Theorems~\ref{thmRHL} and ~\ref{thmRHR}]
We induct on $d$. 
Both statements are trivial when $d = 0$. If Theorems~\ref{thmRHL} and~\ref{thmRHR} hold for all Lefschetz modules of degree at most $d-1$, then Propositions~\ref{prop:raised}, \ref{prop:easyprimitive}, and \ref{prop:criticalprimitive} imply that the conclusions of Theorem~\ref{thmRHL} and Theorem~\ref{thmRHR} hold for all $\eta \in \mathscr{K}_A$. The action of $B^1 \subseteq A^1$ on $\operatorname{Gr}$ is $0$. Then the result follows, as $\eta \in \mathscr{K}_{A/B}$ if and only if there is $b \in B^1$ such that $\eta + b \in \mathscr{K}_A$. 
\end{proof}

\begin{proof}[Proofs of Corollaries~\ref{cor:kernellefschetz} and ~\ref{cor:kernellefschetz2}]
It follows from Corollary~\ref{cor:lefschetzforell} that
$\operatorname{Prim}_\eta^{\bullet,j}$ satisfies the hard Lefschetz property with respect to $\mathscr{K}_B$. 
Let $\mathscr{R}^j$ be the $B$-invariant symmetric bilinear form
\[
\operatorname{Prim}_\eta^{\bullet,j} \times \operatorname{Prim}_\eta^{\bullet,j} \longrightarrow \mathbb{R}, \qquad (x_1,x_2)\longmapsto \underline{\mathscr{Q}}(x_1, \eta^{d-j} * x_2).
\]
It follows from Theorem~\ref{thmRHL} and Proposition~\ref{prop:nondegenerate} that,
for any nonnegative $i \le j/2$ and $\ell \in \mathscr{K}_B$, the symmetric bilinear pairing
\[
\operatorname{Gr}^{i, j} \times \operatorname{Gr}^{i, j} \longrightarrow \mathbb{R}, \qquad 
(x_1,x_2) \mapsto \underline{\mathscr{Q}}(x_1, \eta^{d-j} * \ell^{j - 2i}x_2)
\]
is nondegenerate. Combined with the hard Lefschetz property of $\operatorname{Prim}_\eta^{\bullet,j}$, this implies Poincar\'e duality for $(\operatorname{Prim}_\eta^{\bullet,j},\mathscr{R}^j)$.
The Hodge--Riemann relations for $(\operatorname{Prim}_\eta^{\bullet,j},\mathscr{R}^j)$ over $(B,\mathscr{K}_B)$ follow from Theorem~\ref{thmRHR}.

Similarly, it follows from Theorem~\ref{thmRHL} and Corollary~\ref{cor:lefschetzforell} that $\operatorname{Prim}_\ell^{\bullet,[k]}$ satisfies the hard Lefschetz property with respect to $\mathscr{K}_{A/B}$. Let $\mathscr{R}^{[k]}$ be the $A$-invariant symmetric bilinear pairing
\[
\operatorname{Prim}_\ell^{\bullet,[k]} \times \operatorname{Prim}_\ell^{\bullet,[k]} \longrightarrow \mathbb{R}, \qquad (x_1,x_2)\longmapsto \underline{\mathscr{Q}}(x_1, \ell^k x_2).
\]
Poincar\'{e} duality for $(\operatorname{Prim}_\ell^{\bullet,[k]},\mathscr{R}^{[k]})$ follows from Theorem~\ref{thmRHL} and Proposition~\ref{prop:nondegenerate} as before, and the Hodge--Riemann relations  for $(\operatorname{Prim}_\ell^{\bullet,[k]},\mathscr{R}^{[k]})$ over $(A,\mathscr{K}_{A/B})$ follow from Theorem~\ref{thmRHR}. 
\end{proof}

\section{Proof of the decomposition package}\label{sec:proof}

In this section, we prove Theorem~\ref{thmDecomposition}, from which Theorems~\ref{thmHom},~\ref{thmPD},~\ref{thmHL}, and ~\ref{thmHR} follow. 
Let $(A,\mathscr{K}_A)$ be as in Section~\ref{sec:introduction}. We first discuss some general properties of Lefschetz modules. 
We write $\operatorname{Mod}_A$ for the category of finite-dimensional graded $A$-modules with degree-preserving $A$-module homomorphisms. 
For a nonnegative integer $d$, let $\operatorname{L}^d(\mathscr{K}_A)$  be the full subcategory of $\operatorname{Mod}_A$ consisting of those $M$ for which there is a real-valued $A$-invariant symmetric bilinear form $\mathscr{Q}$ on $M$ 
such that $(M,\mathscr{Q})$ is a Lefschetz module of degree $d$ over ($A,\mathscr{K}_A)$.

\begin{proposition}\label{prop:category}
 $\operatorname{L}^d(\mathscr{K}_A)$ is a semisimple  abelian subcategory of $\operatorname{Mod}_A$.
\end{proposition}

We record the following formal implications of Proposition~\ref{prop:category}:
\begin{enumerate}[(1)]\itemsep 5pt
\item Any graded $A$-module summand of an object $M$ in $\operatorname{L}^d(\mathscr{K}_A)$ is an object of $\operatorname{L}^d(\mathscr{K}_A)$, as it is the kernel of a graded $A$-module homomorphism $M \to M$.
\item The kernel of any graded $A$-module homomorphism $M \to N$ with $M$ and $N$ in $\operatorname{L}^d(\mathscr{K}_A)$ is a graded $A$-module summand of $M$. 
\item An object of $\operatorname{L}^d(\mathscr{K}_A)$ is simple in $\operatorname{L}^d(\mathscr{K}_A)$ 
if and only if it is indecomposable in $\operatorname{Mod}_A$. 
\item Any graded $A$-module homomorphism between two indecomposable graded $A$-modules in $\operatorname{L}^d(\mathscr{K}_A)$ is either zero or an isomorphism.
\end{enumerate}

For an object $M$ of $\operatorname{Mod}_A$, the \emph{graded dual} $M^{\vee}$ is the module $\bigoplus_{k \in \mathbb{Z}} (M^k)^*[k]$, where $(M^k)^*$ is the dual vector space to $M^k$. This is a graded $A$-module, with module structure given by $(a \cdot f)(m) = f(a \cdot m)$.

\begin{proof}
We first show that  $\operatorname{L}^d(\mathscr{K}_A)$ is an abelian subcategory of $\operatorname{Mod}_A$. Clearly, if $(M_1,\mathscr{Q}_1)$ and $(M_2,\mathscr{Q}_2)$ are Lefschetz modules of degree $d$ over $(A,\mathscr{K}_A)$, then $(M_1 \oplus M_2,\mathscr{Q}_1 \oplus \mathscr{Q}_2)$ is a  Lefschetz module of degree $d$ over $(A,\mathscr{K}_A)$. 

We check that, if $M$ and $N$ are in $\operatorname{L}^d(\mathscr{K}_A)$, then the kernel of any graded $A$-module homomorphism $\varphi \colon M \to N$ is in $\operatorname{L}^d(\mathscr{K}_A)$. 
Let $K$ be the kernel of $\varphi$, and choose $\mathscr{Q}$ so that $(M,\mathscr{Q})$ is a Lefschetz module of degree $d$ over $(A,\mathscr{K}_A)$. We claim that $(K,\mathscr{Q})$ is a Lefschetz module of degree $d$ over $(A,\mathscr{K}_A)$.
 For this, it is enough to show that the restriction of $\mathscr{Q}$ to $K$ is nondegenerate. Indeed, this implies that $\dim K^i=\dim K^{d-i}$ for all $i$, and hence the hard Lefschetz property and the Hodge--Riemann relations for $(M,\mathscr{Q})$ imply those for $(K,\mathscr{Q})$. 
 Let $x$ be any nonzero element of $K^i$, let $c$ be the minimum of $i$ and $d-i$, and let $\eta$ be any element of $\mathscr{K}_A$. By the hard Lefschetz property of $M$ with respect to $\eta$, we have the Lefschetz decomposition
 \[
x=\eta^i y_0+\eta^{i-1} y_1 +\cdots+ \eta^{i-c} y_{c}, 
 \]
 where $y_j$ is an element of $M^j$ which is annihilated by $\eta^{d-2j+1}$. 
 This gives
 \[
0=\eta^i \varphi(y_0)+\eta^{i-1} \varphi(y_1) +\cdots+ \eta^{i-c} \varphi(y_{c}),  
 \]
and the uniqueness of the Lefschetz decomposition of $N$ with respect to $\eta$ shows that $y_j \in K^j$ for all $j$. By the Hodge--Riemann relations for $(M,\mathscr{Q})$ and the orthogonality of the primitive Lefschetz summands of $x$ with respect to $\eta$, for any nonzero $y_j$, we have
\[
\mathscr{Q}(x,\eta^{d-i-j} y_j)=\sum_{k=0}^c \mathscr{Q}(\eta^{i-k} y_k ,\eta^{d-i-j} y_j)=\mathscr{Q}(\eta^{i-j} y_j ,\eta^{d-i-j} y_j)=\mathscr{Q}(y_j ,\eta^{d-2j} y_j) \neq 0.
\]
Since $\eta^{d-i-j} y_j$ is an element of $K^{d-i}$, this shows that the restriction of $\mathscr{Q}$ to $K$ is nondegenerate.

We check that, if $M$ and $N$ are in $\operatorname{L}^d(\mathscr{K}_A)$, then the cokernel of any graded $A$-module homomorphism $\varphi \colon M \to N$ is in $\operatorname{L}^d(\mathscr{K}_A)$.
Choose $\mathscr{Q}$ and $\mathscr{R}$ so that $(M,\mathscr{Q})$ and $(N,\mathscr{R})$ are Lefschetz modules of degree $d$ over $(A,\mathscr{K}_A)$.
Using  Poincar\'e duality for $(M,\mathscr{Q})$ and $(N,\mathscr{R})$, we may define $\psi \colon N \to M$ by the property that
\[
\mathscr{Q}(x,\psi(y))=
\mathscr{R}(\varphi(x),y) \ \ \text{for all $x \in M$ and $y \in N$.}
\]
Since $\mathscr{Q}$ and $\mathscr{R}$  are $A$-invariant, $\psi$ is a graded $A$-module homomorphism. 
By the previous argument, the kernel of $\psi$ is in  $\operatorname{L}^d(\mathscr{K}_A)$. 
Since $\ker(\psi)^\vee[-d] \simeq \ker(\psi)$ in $\operatorname{Mod}_A$ by Poincar\'e duality,
$\operatorname{coker}(\varphi)\simeq \ker(\psi)^\vee[-d]$
 is in  $\operatorname{L}^d(\mathscr{K}_A)$.

We next show that $\operatorname{L}^d(\mathscr{K}_A)$ is semisimple. It is enough to show that the kernel of any graded $A$-module homomorphism $\varphi \colon M \to N$ is a direct summand of $M$ for any $M$ and $N$ in $\operatorname{L}^d(\mathscr{K}_A)$. 
 Choose $\mathscr{Q}$ so that $(M,\mathscr{Q})$ is a Lefschetz module of degree $d$ over $(A,\mathscr{K}_A)$.  Let $K$ be the kernel of $\varphi$, and let $K^\perp$ be the orthogonal complement of $K$ in $M$ with respect to $\mathscr{Q}$. Since $\mathscr{Q}$ is $A$-invariant, $K^\perp$ is a graded $A$-submodule of $M$. 
By the previous argument, the restriction of $\mathscr{Q}$ to $K$ is nondegenerate, so $M = K \oplus K^\perp$. 
\end{proof}

\begin{corollary}\label{cor:hom}
If $M$ is an indecomposable graded $A$-module in $\operatorname{L}^d(\mathscr{K}_A)$, then 
there is an isomorphism of $\mathbb{R}$-algebras
\[
\operatorname{Hom}_A(M,M)\simeq \mathbb{R}, \mathbb{C}, \ \text{or} \ \mathbb{H}.
\]
If $M$ and $N$ are nonisomorphic indecomposable graded $A$-modules in $\operatorname{L}^d(\mathscr{K}_A)$, then
\[
\operatorname{Hom}_A(M,N)=0. 
\]
\end{corollary}

Corollary~\ref{cor:hom} implies that, for example, an indecomposable graded $A$-module in $\operatorname{L}^{2d}(\mathscr{K}_A)$ is either $1$-dimensional or has no nonzero socle in degree  $d$. 

\begin{proof}[Proof of Corollary~\ref{cor:hom}]
This follows from Schur's lemma applied to $\operatorname{L}^d(\mathscr{K}_A)$, combined with Frobenius' theorem on finite-dimensional division algebras over $\mathbb{R}$.
\end{proof}

For any $M$ in $\operatorname{Mod}_A$, consider the vector space $\mathscr{S}^d_M$ of all
real-valued $A$-invariant symmetric bilinear forms $\mathscr{Q}$ on $M$ satisfying the degree orthogonality relations 
\[
\mathscr{Q}(M^i,M^j)=0 \ \ \text{whenever $i+j\neq d$}.
\]
Let $\mathscr{C}^d_M=\mathscr{C}^d_M(\mathscr{K}_A)$ be the subset of all $\mathscr{Q}$  in $\mathscr{S}^d_M$ such that $(M,\mathscr{Q})$ is a Lefschetz module of degree $d$ over $(A,\mathscr{K}_A)$. 
Note that the Hodge--Riemann relations impose a convex condition on $\mathscr{Q}$: If $(M,\mathscr{Q}_1)$ and $(M,\mathscr{Q}_2)$ satisfy the Hodge--Riemann relations with respect to $\mathscr{K}_A$, then so does $(M,\mathscr{Q}_1+\mathscr{Q}_2)$. 
This implies that 
$\mathscr{C}^d_M$ is a convex cone in $\mathscr{S}^d_M$.

\begin{proposition}\label{prop:uniqueform}
 If $M$ is an indecomposable graded $A$-module in $\operatorname{L}^d(\mathscr{K}_A)$, then $\dim \mathscr{S}^d_M=1$.
\end{proposition}

\begin{proof}
Note that, assuming  Poincar\'e duality and the hard Lefschetz property for $(M,\mathscr{Q})$, 
the validity of the Hodge--Riemann relations for $(M,\mathscr{Q})$ is equivalent to the statement that, for any nonnegative $k \le \frac{d}{2}$, $\eta \in \mathscr{K}_A$, and $\mathscr{R} \in \mathscr{C}^d_M$, we have
\[
\text{signature}\Big(\mathscr{Q}(-,\eta^{d-2k}-) \colon M^k \times M^k \to \mathbb{R}\Big)=
\text{signature}\Big(\mathscr{R}(-,\eta^{d-2k}-) \colon M^k \times M^k \to \mathbb{R}\Big).
\]
The hard Lefschetz property depends only on the graded $A$-module structure of $M$. Therefore a limit point $\mathscr{Q}$ of $\mathscr{C}^d_M$ lies in $\mathscr{C}^d_M$ if and only if $\mathscr{Q}$ satisfies Poincar\'{e} duality. 
This implies that $\mathscr{C}^d_M$ is open. 

Suppose that $\mathscr{S}^d_M$ has dimension larger than $1$.
In this case, since $\mathscr{C}^d_M \cap -\mathscr{C}^d_M =\varnothing$,
there is a nonzero $\mathscr{Q}$ in the closure of $\mathscr{C}^d_M$ but not in $\mathscr{C}^d_M$.\footnote{This is the step which fails in Examples~\ref{ex:Q-nonuniqueness} and ~\ref{ex:Q-nonuniqueness2} over $\mathbb{Q}$. In those cases, all the rays in the boundary of the cone $\mathscr{C}^d_M$ are irrational.} 
This $\mathscr{Q}$ must fail Poincar\'e duality, and it  follows that the graded $A$-module homomorphism
\[
\varphi \colon M \to M^\vee[-d], \quad x \longmapsto \big(y \longmapsto \mathscr{Q}(x,y)\big)
\]
is nonzero and has a nonzero kernel. Since $M$ is in $\operatorname{L}^d(\mathscr{K}_A)$, there is an isomorphism of graded $A$-modules  $M^\vee[-d] \simeq M$. It follows from Proposition~\ref{prop:category} that the kernel of $\varphi$ is a nonzero proper graded $A$-module summand of $M$, contradicting the indecomposability of $M$.
\end{proof}

We conclude the preliminary material with two simple facts relating $\operatorname{L}^d(\mathscr{K}_A)$ for different values of $d$.

\begin{proposition}\label{prop:shift}
If $M\in \operatorname{L}^d(\mathscr{K}_A)$ and  $M[-k]$ is nonnegatively graded, then
$M[-k] \in \operatorname{L}^{d+2k}(\mathscr{K}_A)$. 
\end{proposition}

\begin{proof}
If $(M,\mathscr{Q})$ is a Lefschetz module of degree $d$ over $(A,\mathscr{K}_A)$, then $(-1)^k\mathscr{Q}[-k]$ gives $M[-k]$ the structure of a Lefschetz module of degree $d+2k$ over $(A,\mathscr{K}_A)$.
\end{proof}

\begin{proposition}\label{prop:van}
If $M$ is in $\operatorname{L}^d(\mathscr{K}_A)$ and $N$ is in $\operatorname{L}^e(\mathscr{K}_A)$ with $d<e$, then
\[
\operatorname{Hom}_A(M,N)=0.
\]
\end{proposition}

\begin{proof}
Let $\varphi \colon M \to N$ be a graded $A$-module homomorphism.
Let $x$ be any nonzero element of $M^i$, let $c$ be the minimum of $i$ and $d-i$, and let $\eta$ be any element of $\mathscr{K}_A$. By the hard Lefschetz property of $M$ with respect to $\eta$, we have the Lefschetz decomposition
 \[
x=\eta^i y_0+\eta^{i-1} y_1 +\cdots+ \eta^{i-c} y_{c}, 
 \]
 where $y_j$ is an element of $M^j$ annihilated by $\eta^{d-2j+1}$. 
 Thus $\varphi(y_j)$ is an element of $N^j$ annihilated by $\eta^{d-2j+1}$, and hence by $\eta^{e-2j}$. The hard Lefschetz property of $N$ implies that $\varphi(y_j)=0$ for all $j$, so we conclude that $\varphi(x)=0$.
\end{proof}

Let $(B,\mathscr{K}_B)$  and $(C,\mathscr{K}_C)$ be as in Section~\ref{sec:introduction}. 
Recall that $C$ is a graded $B$-subalgebra of $A$ that preserves the perverse filtration of $M$ over $B$, so $\operatorname{Gr}$ has the structure of a graded $C$-module such that $C^k \operatorname{Gr}^{i,j} \subseteq \operatorname{Gr}^{i + k, j}$.

\begin{proposition}\label{prop:decomposition}
There is an isomorphism of  graded $C$-modules $M\simeq \operatorname{Gr}$. 
\end{proposition}

The proof follows an argument of Deligne \cite{Deligen1994} and produces a distinguished isomorphism for each $\eta \in \mathscr{K}_{A/B}$. 
See \cite{MR2063103,dCCanonical} for related discussions.

\begin{proof}
Choose $\eta \in \mathscr{K}_{A/B}$, and set $M_0 \coloneq M$. Using  Theorem~\ref{thmRHL}, we can produce a splitting of the inclusion $P_0 \hookrightarrow M_0$ by the composition
\begin{center}
\begin{tikzcd}
M_0 \arrow[r, "\eta^d"] & M_0 \arrow[r] & \operatorname{Gr}^{\bullet, 2d} \arrow[r, "(\eta^d)^{-1}"] & \operatorname{Gr}^{\bullet, 0} =P_0.
\end{tikzcd}
\end{center}
This gives a decomposition of graded $C$-modules 
\[
M_0 = P_0 \oplus M_{0,1} = \operatorname{Gr}^{\bullet, 0} \oplus M_{0,1},
\]
where $M_{0,1} \simeq M_0/P_0$ is the kernel of the displayed composition. 
The perverse filtration on $M_0$ induces a filtration on $M_{0,1}$ whose associated graded is $\oplus_{j = 1}^{2d} \operatorname{Gr}^{\bullet, j}$. 
 Using  Theorem~\ref{thmRHL} again, we can produce a splitting of $M_{0,1} \to M_0 \to \operatorname{Gr}^{\bullet, 2d}$ by the composition
\begin{center}
\begin{tikzcd}
\operatorname{Gr}^{\bullet, 2d} \arrow[r, "(\eta^d)^{-1}"] & \operatorname{Gr}^{\bullet, 0} \arrow[r] & M_0 \arrow[r, "\eta^d"] & M_0 \arrow[r] & M_{0,1}.
\end{tikzcd}
\end{center}
This gives a decomposition of graded $C$-modules $M_{0,1} \simeq \operatorname{Gr}^{\bullet, 2d} \oplus M_1$, 
and hence
\[
M_0\simeq \operatorname{Gr}^{\bullet, 0} \oplus M_{0,1}\simeq \operatorname{Gr}^{\bullet, 0} \oplus \operatorname{Gr}^{\bullet, 2d} \oplus M_1, \ \ \text{where $M_1\simeq P_{2d-1}/P_0$.}
\]
The filtration on $M_{0,1}$ induces a filtration on $M_1$ whose associated graded is 
$\oplus_{j = 1}^{2d-1} \operatorname{Gr}^{\bullet, j}$. 

 Suppose that, after
$r$ steps, we have constructed a $C$-module summand $M_r$ of 
$M$ with a filtration whose associated graded is $\oplus_{j = r}^{2d-r} \operatorname{Gr}^{\bullet, j}$. We use Theorem~\ref{thmRHL} to split $\operatorname{Gr}^{\bullet, r}\to M_r
$ by
\[
\begin{tikzcd}
M_r \arrow[r, "\eta^{d-r}"] 
& M \arrow[r] 
& M_r \arrow[r] 
& \operatorname{Gr}^{\bullet, 2d-r} \arrow[r, "(\eta^{d-r})^{-1}"] 
& \operatorname{Gr}^{\bullet, r}.
\end{tikzcd}
\]
This gives a graded $C$-module 
decomposition 
\[
M_r=\operatorname{Gr}^{\bullet, r} \oplus M_{r,r+1},
\]
and the filtration on $M_r$ induces  a filtration on $M_{r,r+1}$ 
with associated graded $\oplus_{j = r+1}^{2d-r} \operatorname{Gr}^{\bullet, j}$. 
 Next we use Theorem~\ref{thmRHL} to split
 $M_{r,r+1} \to  \operatorname{Gr}^{\bullet, 2d-r}$ by
\[
\begin{tikzcd}
\operatorname{Gr}^{\bullet, 2d-r} \arrow[r, "(\eta^{d-r})^{-1}"] & \operatorname{Gr}^{\bullet, r} \arrow[r] & M_r \arrow[r, "\eta^{d-r}"]&M \ar[r] & M_r \arrow[r] & M_{r,r+1}.
\end{tikzcd}
\]
This gives a graded $C$-module decomposition 
\[
M_{r,r+1}\simeq  \operatorname{Gr}^{\bullet,2d-r} \oplus M_{r+1},
\]
and the filtration on $M_{r,r+1}$ induces 
 a filtration on $M_{r+1}$ with associated graded  $\oplus_{j = r+1}^{2d-r-1} \operatorname{Gr}^{\bullet, j}$.

After the step $r=d-1$, the remaining summand is isomorphic to $\operatorname{Gr}^{\bullet,d}$. Combining these decompositions gives a graded $C$-module isomorphism $M \simeq \operatorname{Gr}$.
\end{proof}

Proposition~\ref{prop:decomposition} and the decomposition of $M$ in \eqref{eq:decomposition-indecomposable2} give a graded $C$-module isomorphism
\begin{equation}\label{eq:Gr}
\operatorname{Gr}\simeq \bigoplus_{\alpha} \bigoplus_k L_\alpha[-k]^{\oplus \mu(\alpha,k)}.
\end{equation}
We will compare this decomposition with the Lefschetz decomposition of $\operatorname{Gr}$  obtained from Theorem~\ref{thmRHL}: 
For any $\eta \in \mathscr{K}_{A/B}$, we have an identification of graded $C$-modules
\begin{equation}\label{eq:prim}
\operatorname{Gr}=\bigoplus_{0 \le j \le d} \bigoplus_{0 \le i \le d-j}  \eta^{i} * \operatorname{Prim}^{\bullet,j}_{\eta}, 
\end{equation}
where $\operatorname{Prim}^{\bullet,j}_{\eta}$ is the kernel appearing in Corollary~\ref{cor:kernellefschetz}.

\begin{proof}[Proof of Theorem~\ref{thmDecomposition}]
We first check that the perverse filtration on $M$ over $B$ is the same as the perverse filtration on $M$ over the subalgebra of $C$ generated by $C^1$. The perverse filtration on $M$ over the subalgebra generated by $C^1$ is the unique $C^1$-stable filtration for which some element of $C^1$ satisfies the hard Lefschetz property on the associated graded. The perverse filtration on $M$ over $B$ has this property. 

By replacing $B$ with the subalgebra of $C$ generated by $C^1$, we may assume that $B^1=C^1$ and  $\mathscr{K}_B = \mathscr{K}_C$. Under this assumption, 
Corollary~\ref{cor:kernellefschetz} states that 
$\operatorname{Prim}^{\bullet,j}_\eta$ is in $\operatorname{L}^j(\mathscr{K}_C)$. 

Comparing \eqref{eq:Gr} and \eqref{eq:prim}, 
 we deduce from the Krull--Schmidt theorem for $\operatorname{Mod}_C$ that every indecomposable graded $C$-module summand of $M$ must appear as a summand of $\eta^i* \operatorname{Prim}_{\eta}^{\bullet,j}$ for some $i$ and $j$. 
 By Theorem~\ref{thmRHL}, we have
 \[
 \eta^i* \operatorname{Prim}_{\eta}^{\bullet,j} \simeq \operatorname{Prim}_{\eta}^{\bullet,j}[-i] \ \ \text{in $\operatorname{Mod}_C$.}
 \]
Therefore, for each $\alpha$, any summand $L_\alpha[-k]$ of $M$ with the smallest possible $k$ among those with nonzero $\mu(\alpha,k)$ must be a summand of $\operatorname{Prim}_{\eta}^{\bullet,j}$ for some $j$. 
Since 
$\operatorname{Prim}^{\bullet,j}_\eta$ is in $\operatorname{L}^j(\mathscr{K}_C)$ by Corollary~\ref{cor:kernellefschetz},
 it follows from Proposition~\ref{prop:category} that $L_\alpha[-k]$ is in $\operatorname{L}^j(\mathscr{K}_C)$.

By Proposition~\ref{prop:shift}, this means that $L_\alpha$ is in $\operatorname{L}^{e(\alpha)}(\mathscr{K}_C)$.
In other words, there is a $C$-invariant symmetric bilinear form $\mathscr{Q}_\alpha$ on $L_\alpha$ satisfying the degree orthogonality relations such that $(L_\alpha,\mathscr{Q}_\alpha)$ is a Lefschetz module of degree $e(\alpha)$ over $(C,\mathscr{K}_C)$. 
Since $L_\alpha$ is indecomposable in $\operatorname{Mod}_C$,  Proposition~\ref{prop:uniqueform} applies, and this $\mathscr{Q}_{\alpha}$ is the unique $C$-invariant symmetric bilinear form on $L_\alpha$ satisfying the degree orthogonality relations, up to a nonzero real multiple.
By Corollary~\ref{cor:hom}, there is an isomorphism of $\mathbb{R}$-algebras
\[
\operatorname{Hom}_C(L_\alpha,L_\alpha) \simeq \mathbb{R}, \mathbb{C}, \ \text{or} \ \mathbb{H}.
\]
By the same token, 
if $L_\alpha$ and $L_\beta$ are not isomorphic in $\operatorname{Mod}_C$ and $e(\alpha)=e(\beta)+2k$, then  
\[
\operatorname{Hom}_C(L_\alpha,L_\beta[-k])=0.
\]
By Proposition~\ref{prop:van}, if $\alpha=\beta$ and $k$ is  positive, or if $\alpha\neq \beta$ and $e(\alpha)<e(\beta)+2k$, then
\[
\operatorname{Hom}_C(L_\alpha,L_\beta[-k])=0. 
\]
This completes the proof.
\end{proof}

Theorems~\ref{thmHom},~\ref{thmPD},~\ref{thmHL}, and ~\ref{thmHR} follow from the case $C=B$.

\bibliography{refs}
\bibliographystyle{amsalpha}

\end{document}